\documentclass[11pt]{article}

\usepackage{amscd,amsmath, amssymb, fancyhdr, epsfig,url,color}

\numberwithin{equation}{section}

\newcommand{\version}{version 2.3.5,\ \  May 11, 2021}

\def\eqref#1{(\ref{#1})}
\newcommand{\goth}{\mathfrak}

\newcommand{\arrow}{{\:\longrightarrow\:}}
\newcommand{\Z}{{\mathbb Z}}
\newcommand{\C}{{\mathbb C}}

\newcommand{\R}{{\mathbb R}}
\newcommand{\Q}{{\mathbb Q}}

\def\1{\sqrt{-1}\:}

\newcommand{\cntrct}                
{\hspace{2pt}\raisebox{1pt}{\text{$\lrcorner$}}\hspace{2pt}}

\makeatletter
\def\x@arrow{\DOTSB\Relbar}
\def\xlongequalsignfill@{\arrowfill@\x@arrow\Relbar\x@arrow}
\newcommand{\xlongequal}[2][]{%
        \ext@arrow 0099\xlongequalsignfill@{#1}{#2}}
\def\xlongrightarrowfill@{\arrowfill@\relbar\relbar\longrightarrow}
\newcommand{\xlongrightarrow}[2][]{%
        \ext@arrow 0099\xlongrightarrowfill@{#1}{#2}}
\makeatother

\renewcommand{\bar}{\overline}
\renewcommand{\phi}{\varphi}
\renewcommand{\epsilon}{\varepsilon}
\renewcommand{\geq}{\geqslant}
\renewcommand{\leq}{\leqslant}
\renewcommand{\max}{{\rm max}}

\newcommand{\Vol}{\operatorname{Vol}}

\newcommand{\Lie}{\operatorname{Lie}}
\newcommand{\Sp}{\operatorname{Sp}}
\newcommand{\Diff}{\operatorname{Diff}}

\newcommand{\Teich}{\operatorname{\sf Teich}}

\newcommand{\Comp}{\operatorname{\sf Comp}}

\newcommand{\Per}{\operatorname{\sf Per}}
\newcommand{\Kah}{\operatorname{\sf Kah}}

\newcounter{Mycounter}[section]
\newcounter{lemma}[section]
\renewcommand{\thelemma}{{Lemma \thesection.\arabic{lemma}}}
\newcommand{\lemma}{%
    \setcounter{lemma}{\value{Mycounter}}
    \refstepcounter{lemma}
    \stepcounter{Mycounter}
    {\noindent \bf \thelemma:\ }}

\newcounter{claim}[section]
\renewcommand{\theclaim}{{Claim \thesection.\arabic{claim}}}
\newcommand{\claim}{%
    \setcounter{claim}{\value{Mycounter}}
    \refstepcounter{claim}
    \stepcounter{Mycounter}
    {\noindent \bf \theclaim:\ }}

\newcounter{sublemma}[section]
\newcounter{corollary}[section]
\renewcommand{\thecorollary}{{Corollary \thesection.\arabic{corollary}}}
\newcommand{\corollary}{%
    \setcounter{corollary}{\value{Mycounter}}
    \refstepcounter{corollary}
    \stepcounter{Mycounter}
    {\noindent \bf \thecorollary:\ }}

\newcounter{theorem}[section]
\renewcommand{\thetheorem}{{Theorem \thesection.\arabic{theorem}}}
\newcommand{\theorem}{%
    \setcounter{theorem}{\value{Mycounter}}
    \refstepcounter{theorem}
    \stepcounter{Mycounter}
    {\noindent \bf \thetheorem:\ }}

\newcounter{conjecture}[section]
\renewcommand{\theconjecture}{{Conjecture \thesection.\arabic{conjecture}}}
\newcommand{\conjecture}{%
    \setcounter{conjecture}{\value{Mycounter}}
    \refstepcounter{conjecture}
    \stepcounter{Mycounter}
    {\noindent \bf \theconjecture:\ }}

\newcounter{proposition}[section]
\renewcommand{\theproposition}
      {{Proposition \thesection.\arabic{proposition}}}
\newcommand{\proposition}{%
    \setcounter{proposition}{\value{Mycounter}}
    \refstepcounter{proposition}
    \stepcounter{Mycounter}
    {\noindent \bf \theproposition:\ }}

\newcounter{definition}[section]
\renewcommand{\thedefinition}
      {{Definition~\thesection.\arabic{definition}}}
\newcommand{\definition}{%
    \setcounter{definition}{\value{Mycounter}}
    \refstepcounter{definition}
    \stepcounter{Mycounter}
    {\noindent \bf \thedefinition:\ }}

\newcounter{example}[section]
\newcounter{remark}[section]
\renewcommand{\theremark}{{Remark \thesection.\arabic{remark}}}
\newcommand{\remark}{%
    \setcounter{remark}{\value{Mycounter}}
    \refstepcounter{remark}
    \stepcounter{Mycounter}
    {\noindent \bf \theremark:\ }}

\newcounter{problem}[section]
\newcounter{question}[section]
\renewcommand{\thequestion}{{Question \thesection.\arabic{question}}}
\newcommand{\question}{%
    \setcounter{question}{\value{Mycounter}}
    \refstepcounter{question}
    \stepcounter{Mycounter}
    {\noindent \bf \thequestion:\ }}

\makeatletter

\@addtoreset{equation}{section} \@addtoreset{footnote}{section} \makeatother

\def\blacksquare{\hbox{\vrule width 5pt height 5pt depth 0pt}}
\def\endproof{\blacksquare}

\newcommand{\T}{{\mathbb{T}}}
\newcommand{\N}{{\mathbb{N}}}

\newcommand{\bJ}{{\overline{J}}}

\newcommand{\bomega}{{\overline{\omega}}}

\newcommand{\bw}{{\bar{w}}}
\newcommand{\bbw}{{\bar{\bf w}}}
\newcommand{\bGamma}{{\overline{\Gamma}}}

\newcommand{\cB}{{\mathcal{B}}}
\newcommand{\cD}{{\mathcal{D}}}
\newcommand{\cC}{{\mathcal{C}}}

\newcommand{\cF}{{\mathcal{F}}}

\newcommand{\cT}{{\mathcal{T}}}
\newcommand{\cU}{{\mathcal{U}}}

\newcommand{\cX}{{\mathcal{X}}}

\newcommand{\hcC}{{\hat{\mathcal{C}}}}

\newcommand{\hcT}{{\hat{\mathcal{T}}}}

\newcommand{\gog}{{\mathfrak{g}}}
\newcommand{\gogc}{{\mathfrak{g}_\C}}
\newcommand{\goh}{{\mathfrak{h}}}
\newcommand{\gohc}{{\mathfrak{h}_\C}}
\newcommand{\gos}{{\mathfrak{s}}}
\newcommand{\gosc}{{\mathfrak{s}_\C}}
\newcommand{\gok}{{\mathfrak{k}}}

\newcommand{\gogl}{{\mathfrak{gl}}}

\newcommand{\fsl}{\mathfrak{sl}}
\newcommand{\fgl}{\mathfrak{gl}}

\newcommand{\Span}{{{\rm Span}}}
\newcommand{\Camp}{{{\rm Camp}\,}}

\newcommand{\tI}{{\widetilde{I}}}
\newcommand{\tJ}{{\widetilde{J}}}

\newcommand{\tM}{{\widetilde{M}}}

\newcommand{\tomega}{{\widetilde{\omega}}}

\newcommand{\teta}{{\widetilde{\eta}}}
\newcommand{\talpha}{{\widetilde{\alpha}}}

\begin{document}

%
\begin{center}
{\LARGE\bf
Unobstructed symplectic packing for tori and hyperk\"ahler manifolds\\[4mm]
}

Michael Entov\footnote{Partially supported by the Israel Science
Foundation grant $\#$ 1096/14 and by M. \& M. Bank Mathematics
Research Fund.}, Misha Verbitsky\footnote{Partially supported
by RSCF grant 14-21-00053 within AG Laboratory NRU-HSE.
 }

\end{center}


\begin{abstract}
\small Let
$M$ be a closed symplectic manifold of volume $V$. We say that
the symplectic packings of $M$ by balls are unobstructed if any collection
of disjoint symplectic balls (of possibly different radii) of total volume
less than $V$ admits a symplectic embedding to $M$.
 In 1994 McDuff and Polterovich proved that
symplectic packings of K\"ahler manifolds by balls can be characterized in
terms of the K\"ahler cones of their blow-ups. When $M$ is a K\"ahler
manifold which is not a union of its proper subvarieties (such a
manifold is called Campana simple) these K\"ahler cones can be
described explicitly using the Demailly and Paun structure theorem.
We prove that for any Campana simple K\"ahler manifold, as well as for any
manifold which is a limit of Campana simple manifolds in a smooth deformation, the symplectic packings by
balls are unobstructed. This is used to show that the symplectic packings by balls
of all even-dimensional tori equipped with K\"ahler symplectic forms
and of all hyperk\"ahler manifolds of maximal holonomy are
unobstructed. This generalizes a
previous result by Latschev-McDuff-Schlenk. We also consider
symplectic packings by other shapes and show, using Ratner's orbit
closure theorem, that any even-dimensional torus equipped with a
K\"ahler form whose cohomology class is not proportional to a
rational one admits a full symplectic packing by any number
of equal polydisks (and, in particular, by any number of equal
cubes).

\end{abstract}


{\small
\tableofcontents
}


\section{Introduction}


The symplectic packing problem is one of the major
problems of symplectic topology that was introduced, along
with the first results on it, in the famous foundational
paper by Gromov \cite{_Gromov_}.
The most extensively studied version of the problem is the
question about symplectic packings of symplectic
manifolds by balls. In \cite{_McD-Polt_} McDuff and Polterovich reduced the
question about such packings of a symplectic
manifold $(M,\omega)$ to a question about
the structure of the symplectic cone in the cohomology of
a blow-up of $M$. In the same paper they
showed that symplectic packings of K\"ahler manifolds by balls
are deeply related to algebraic geometry that
allows sometimes to describe the shape of the K\"ahler cone in the cohomology of a
K\"ahler manifold.

In this paper we use several strong results from complex geometry in
order to prove the flexibility of symplectic packings by balls for
all even-dimensional tori equipped with K\"ahler symplectic forms as
well as for certain hyperk\"ahler manifolds. Namely, we show that
such a packing is possible as long as the natural volume constraint
is satisfied. In this case we will say that the symplectic packings of the symplectic manifold
by balls are {\it unobstructed}. (We use the term ``unobstructed"
when we consider packings by {\it all} collections of balls of {\it arbitrary
relative sizes}
and the term ``a full symplectic packing" when we talk about packings by
a {\it specific} number of balls, or other shapes, of {\it specific}
relative sizes -- e.g. ``a full symplectic packing by $k$ equal balls").

Our unobstructed symplectic packing theorem
extends a previous result of Latschev-McDuff-Schlenk
\cite{_LMcDS_}.

Our strategy is the following. Let $(M,I,\omega)$ be a closed
connected K\"ahler manifold, where $I$ is the complex structure and
$\omega$ is the symplectic form forming the K\"ahler structure. If
the complex structure $I$ is {\bf Campana simple}, meaning that the
union of the positive-dimensional proper complex submanifolds of
$(M,I)$  is of measure zero, then the Demailly-Paun theorem
\cite{_Demailly_Paun_} allows to give a complete description of the
K\"ahler cones of the relevant blow-ups of $M$. Together with
McDuff-Polterovich theorem mentioned above this implies that the symplectic packings
of $(M,\omega)$ by balls are unobstructed.
Then, {\bf and this is the main novelty of our approach, compared to \cite{_LMcDS_}}, we use the
Kodaira-Spencer stability theorem \cite{_Kod-Spen-AnnMath-1960_} to
show that even if the complex structure $I$ on $M$ is not Campana
simple but can only be approximated in a smooth deformation by Campana simple complex
structures, the result about the K\"ahler cones
in the Campana simple case gives enough information about the {\it
symplectic} cones of the blow-ups of $M$ to yield that
the symplectic packings of $(M,\omega)$ by balls are unobstructed in this case as well
(see
Section~\ref{_Sketch_of_the_proof_of_the_main_thm_on_full_packing_Section_}
for a more detailed outline of this argument and
Section~\ref{_full_packing_Campana_simple_Section_} for a complete
proof). We then use methods of complex geometry to show that if $M$,
$\dim_\R M\geq 4$, is a torus (respectively, a hyperk\"ahler manifold of a
certain type) and $\omega$ is a K\"ahler (respectively, a
hyperk\"ahler) form on $M$, then $I$, appearing with
$\omega$ in a K\"ahler (respectively, a hyperk\"ahler)
structure, on $M$, can be indeed approximated in a smooth deformation
by Campana simple complex structures.

In this paper we also study symplectic packings of the tori and
certain hyperk\"ahler manifolds by arbitrary shapes (with the vanishing second real cohomology).
For such a manifold $M$
we prove the following theorem.
Let $\omega_1$, $\omega_2$ be two K\"ahler (respectively, hyperk\"ahler)
forms on $M$ whose cohomology classes are not proportional to
rational ones. Assume that $\int_M \omega_1^n = \int_M \omega_2^n
>0$, where $2n= \dim_\R M$. In the hyperk\"ahler case assume also that
$\omega_1$ and $\omega_2$ lie in the same deformation class of
hyperk\"ahler forms. Then the maximal fraction of the total volume
that can be filled by packing copies of a given shape of an arbitrary relative size into the
symplectic manifolds $(M,\omega_1)$, $(M,\omega_2)$ is the same.
The proof is based on the ideas from \cite{_Verbitsky:ergodic_}, \cite{_Verbitsky:ICM_}.
Its key step, involving an application of Ratner's orbit
closure theorem, is to show the ergodicity of the action of the group of
orientation-preserving diffeomorphisms of $M$ on the space of K\"ahler forms on
$M$ with a fixed positive total volume.

Combining the latter result with the fact that the symplectic packings of
tori by balls are unobstructed, we show that $T^{2n}$, equipped with any K\"ahler
form $\omega$ whose cohomology class is not proportional
to a rational one, admits a full symplectic packing by any
number
of equal $2n$-dimensional polydisks (that is, $2n$-dimensional
direct products of arbitrary symplectic balls), and, in particular,
by any number of equal $2n$-dimensional cubes.

Let us now recall a few preliminaries and present exact statements of our results.

\section{Preliminaries}
\label{_Preliminaries_Section_}

\bigskip
\noindent {\bf Symplectic and complex structures}.
We will view complex structures as tensors, that is, as integrable almost complex structures.

We say that an almost complex structure $J$ and a differential 2-form
$\omega$ on a smooth manifold $M$ are {\bf
compatible} with
each other if $\omega (\cdot, J\cdot)$ is a
$J$-invariant Riemannian metric on $M$.

Any closed differential 2-form compatible with an almost complex
structure is automatically symplectic\footnote{Every symplectic form
admits compatible almost complex structures \cite{_Gromov_} but does
not necessarily admit compatible complex structures. At the same
time an almost complex structure may not be compatible with any
symplectic form.}.

The compatibility between a {\it complex} structure $J$ and a
symplectic form $\omega$ means exactly that $\omega (\cdot, J\cdot)
+ i\omega (\cdot, \cdot)$ is a K\"ahler metric on $M$.

We will call
a symplectic form {\bf K\"ahler}, if it is compatible with {\it some}
complex structure.

We will say that a complex structure is of {\bf K\"ahler type} if
it is compatible with {\it some} symplectic form.

\bigskip
\noindent {\bf Symplectic forms on tori}. Consider a torus $T^{2n} =
\R^{2n}/\Z^{2n}$ and let $\pi: \R^{2n}\to \R^{2n}/\Z^{2n}=T^{2n}$
be the natural projection.

A differential form (respectively, a complex structure) on
$\R^{2n}$ is called {\bf linear}
if it has constant coefficients with respect to the
standard coordinates on $\R^{2n}$,
or, in other words, if it defines a linear exterior form (respectively,
a linear complex structure) on the vector space $\R^{2n}$.
A linear differential form (respectively, a linear complex structure) on $\R^{2n}$
descends under $\pi$ to a differential form (respectively, a complex structure)
on $T^{2n}$. We will call a differential form (respectively, a complex structure) on $T^{2n}$  {\bf linear}
if it can be obtained in this way.

Any linear symplectic form on $T^{2n}$ is compatible with a linear
complex structure (and thus is K\"ahler). Likewise, any linear
complex structure on $T^{2n}$ is compatible with a linear symplectic
form (and thus is a complex structure of K\"ahler type), since the
same holds for linear symplectic forms and linear complex structures
on $\R^{2n}$. In fact, any K\"ahler form on $T^{2n}$ can be mapped
by a symplectomorphism to a linear symplectic form -- see
\ref{_Any_complex_structure_of_Kahler_type_on_torus_linear_Proposition_}.

\bigskip
\noindent {\bf Hyperk\"ahler manifolds}. There are several
equivalent definitions of a hyperk\"ahler manifold.
Since we study hyperk\"ahler manifolds from the symplectic
viewpoint, here is a definition which is close in spirit to
symplectic geometry: A {\bf hyperk\"ahler manifold} is a manifold
equipped with three complex structures $I_1,I_2,I_3$ satisfying the
quaternionic relations and three symplectic forms
$\omega_1,\omega_2,\omega_3$ compatible, respectively, with
$I_1,I_2,I_3$, so that the three Riemannian metrics $\omega_i
(\cdot, I_i\cdot)$, $i=1,2,3$, coincide. Such a collection of
complex structures and symplectic forms on a manifold is called a
{\bf hyperk\"ahler structure} and will be denoted by
$\goh=\{ I_1, I_2,I_3, \omega_1,\omega_2,\omega_3\}$.

We will say that a symplectic form is
{\bf hyperk\"ahler} and a complex structure is of {\bf
hyperk\"ahler type}, if each of them appears in {\it some} hyperk\"ahler
structure. In particular, any hyperk\"ahler symplectic form is
K\"ahler and any complex structure of hyperk\"ahler type is also of
K\"ahler type.

We say that two hyperk\"ahler forms are {\bf hyperk\"ahler deformation equivalent} if they can be connected by a
smooth path of hyperk\"ahler forms.

The real dimension of a manifold admitting a hyperk\"ahler structure has to be divisible by $4$ -- this follows readily from the fact that
the complex structures $I_1, I_2,I_3$ appearing in a
hyperk\"ahler structure on $M$
induce an action of the quaternions on $TM$.

Here is the complete list of currently known closed manifolds
admitting a hyperk\"ahler structure: $T^{4n}$, K3-surfaces, the
deformations of Hilbert schemes of points of K3-surfaces,
deformations of generalized Kummer manifolds,\footnote{A generalized
Kummer manifold is a fiber of the Albanese map from Hilbert scheme
of a torus $T^4$ to $T^4$.} two more ``sporadic" manifolds due to
O'Grady \cite{_OGrady1_, _OGrady2_}, and finite quotients of direct
products of the examples above.

It follows from the famous Calabi-Yau theorem \cite{_Calabi_,_Yau_} that
any K\"ahler form on a hyperk\"ahler manifold is cohomologous to a
unique hyperk\"ahler form \cite[Theorem
  23.5]{_Huybrechts:lec_}. It was conjectured
that any symplectic form on a hyperk\"ahler manifold is
hyperk\"ahler but even for K3-surfaces this is unknown
\cite{_Donaldson:ellipt_}.

A hyperk\"ahler manifold $(M,\goh)$ is called {\bf irreducible
holomorphically symplectic (IHS)} if $\pi_1 (M)=0$ and $\dim_\C
H^{2,0}_I (M,\C)=1$, where $I$ is any of the three complex
structures appearing in $\goh$ and $H^{2,0}_I (M,\C)$ is the
$(2,0)$-part in the Hodge decomposition of $H^2 (M,\C)$ defined by
$I$ (see Section~\ref{_Hodge_structures_Section_}; for all three
complex structures in $\goh$ the space $H^{2,0}_\cdot = (M,\C)$ has
the same dimension). K3-surfaces, as well as the Hilbert schemes of
points for $T^4$ and for K3-surfaces, are IHS. Any closed
hyperk\"ahler manifold admits a finite covering which is the product
of a torus and several IHS hyperk\"ahler manifolds
\cite{_Bogomolov:decompo_}. The IHS hyperk\"ahler manifolds are also
called {\bf hyperk\"ahler manifolds of maximal holonomy}, because
the holonomy group of a hyperk\"ahler manifold is $\Sp(n)$ (the
group of invertible quaternionic $n\times n$-matrices) -- and not
its proper subgroup -- if and only if it is IHS
\cite{_Besse:Einst_Manifo_}.

\hfill

\section{Main results}
\label{_Main_results_Section_}

\subsection{Unobstructed symplectic packings by balls}

By $\Vol$ we will always denote the symplectic volume of a symplectic manifold.

Let $(M,\omega)$, $\dim_\R M = 2n$, be a closed connected symplectic manifold.
We say that the symplectic packings of $(M,\omega)$ by balls are {\bf unobstructed}, if any finite collection of pairwise disjoint
closed round balls in the standard symplectic $\R^{2n}$ of total volume less than
$\Vol (M,\omega)$ has an open neighborhood that can be symplectically embedded into $(M,\omega)$.

\hfill


\theorem\label{_FSP_main_Theorem_}\\ Let $M$ be a torus $T^{2n}$
with a K\"ahler form $\omega$, or an IHS hyperk\"ahler manifold with
a hyperk\"ahler symplectic form $\omega$. Then the symplectic packings of $(M, \omega)$ by balls are unobstructed.


\hfill


This extends a previous result of Latschev-McDuff-Schlenk \cite{_LMcDS_}.
For the proof see Section~\ref{pf-of-fsp-main-thm_Section}.

\subsection{Symplectic packing by arbitrary shapes}
\label{_arbi_shapes_Subsection_}

Let $(U,\eta)$, $\dim_\R U = 2n$, be an open, possibly disconnected,
symplectic manifold, and let $V\subset U$, $\dim_\R V = 2n$, be a
compact, possibly disconnected, submanifold of $U$ with piecewise
smooth boundary.

By a {\bf symplectic embedding of $(V, \eta)$ in $(M,\omega)$} we mean a symplectic
embedding of an open neighborhood of
$V$ in $(U,\eta)$ to $(M,\omega)$.

Set
$$\nu (M,\omega, V) := \frac{\sup_\alpha \Vol (V,\alpha\eta)}{\Vol (M,\omega)},$$
where the supremum is taken over all $\alpha$
such that $(V, \alpha\eta)$ admits a symplectic embedding into $(M,\omega)$. If there is no such $\alpha$, set
$\nu (M,\omega, V) := -\infty$.

We say that $(M,\omega)$ {\bf can be fully packed by $k$ equal copies of $(V,\eta)$} if $$\nu (M,\omega, W) = 1,$$ where
$W$ is a disjoint union of $k$ equal copies of $(V,\eta)$, or, in other words,
if a disjoint union of $k$ copies of $(V,\alpha\eta)$, $\alpha >0$, admits a symplectic embedding into $(M,\omega)$
if and only if $k\Vol (V,\alpha\eta)<\Vol (M,\omega)$.


\hfill


\theorem\label{_Packing_by_arb_shapes_main_Theorem_}\\
With $V\subset (U,\eta)$ as above, assume that $H^2 (V,\R) = 0$. Let
$M$, $\dim_\R M =: 2n\geq 4$, be either an oriented torus $T^{2n}$ or, respectively, a closed connected oriented manifold admitting IHS hyperk\"ahler structures (compatible
with the orientation). Let $\omega_1$, $\omega_2$ be either
K\"ahler forms on $T^{2n}$ or, respectively, hyperk\"ahler forms on $M$.
Assume that
$\int_M \omega_1^n = \int_M \omega_2^n >0$ and that the
cohomology classes $[\omega_1]$, $[\omega_2]$ are not proportional to rational ones.
In the hyperk\"ahler case assume also that $\omega_1$, $\omega_2$ are hyperk\"ahler deformation equivalent.

Then $\nu (M,\omega_1, V) =
\nu (M,\omega_2, V)$.


\hfill


Equivalently, the claim of the theorem can be stated as follows: $(V,\eta)$ can be symplectically embedded
in $(M,\omega_1)$ if and only if it can be symplectically embedded in $(M,\omega_2)$.

The theorem also immediately implies that
$(M,\omega_1)$ can be fully packed by $k$ equal copies of
$(V,\eta)$ if and only if so can $(M,\omega_2)$.

For the proof of the theorem see Section~\ref{_Packing_by_arb_shapes_main_Section_}.

The strategy of the proof is as follows. Without loss of generality,
we may assume that $\Vol (M,\omega)=1$. The group $\Diff^+$ of
orientation-preserving diffeomorphisms of $M$ acts on the space
$\cF$ of K\"ahler  (respectively, in the hyperk\"ahler case, hyperk\"ahler) forms on $M$ of
total volume 1. The function $\omega\mapsto \nu (M,\omega,
V)$ is clearly invariant under the action. We will show
that this function is lower semicontinuous (with respect to the
$C^\infty$-topology on $\cF$) and that the orbit of $\omega$ under
action of $\Diff^+$ is dense in $\cF$ for $M=T^{2n}$ (respectively, in a connected component $\cF^0$ of $\cF$
containing $\omega$ in the hyperk\"ahler case) as long as the cohomology class $[\omega]$ is not proportional to a rational one. Then, since the orbits
of both $\omega_1$ and $\omega_2$ are dense in $\cF$ (respectively, in $\cF^0$), we get, by the
lower semicontinuity, that $\nu (M,\omega_1, V) \leq \nu
(M,\omega_2, V)$ and $\nu (M,\omega_1, V) \geq
\nu (M,\omega_2, V)$, which means that $\nu (M,\omega_1,
V) = \nu (M,\omega_2, V)$.

Clearly, $\nu$ can be replaced in the theorem by any symplectic invariant of a symplectic manifold
which depends (lower or upper) semicontinuously on
 the symplectic form -- for instance,
by any symplectic capacity defined by means of symplectic embeddings (like the Gromov width). {\bf It would be interesting to find other kinds of semicontinuous symplectic
 invariants}.

As an application of \ref{_Packing_by_arb_shapes_main_Theorem_}, consider the case where
$V$ is the union of $k$ disjoint
translated copies of a polydisk
$$B^{2n_1} (R_1)\times\ldots\times B^{2n_l} (R_l)\subset \R^{2n}$$ for some $n_1+\ldots+n_l = n$ and $R_1,\ldots, R_l>0$.
 (Here $B^{2n_i} (R_i)$, $R_i>0$, is a closed round ball in $\R^{2n_i}$ equipped with the standard symplectic form
$\Omega_{2n_i}$, $i=1,\ldots,l$, and $\Omega_{2n_1}\oplus\ldots \oplus\Omega_{2n_l} = dp\wedge dq$. Accordingly $U$ can be taken to be the standard symplectic $\R^{2n}$).
\ref{_Packing_by_arb_shapes_main_Theorem_} allows to prove the following corollary (the proof
also uses \ref{_FSP_main_Theorem_}).


\hfill


\corollary
\label{_full_packing_of_tori_by_polydisks_Corollary_}
\\
Let $\omega$ be a K\"ahler form on $T^{2n}$, $n\geq 2$, and assume
that the cohomology class $[\omega]$ is not proportional to a rational one. Then for any
$k\in\N$ the symplectic manifold $(T^{2n},\omega)$ can be fully
packed by
$k$ equal polydisks $B^{2n_1} (R_1)\times\ldots\times B^{2n_l}
(R_l)$ (for any $n_1,\ldots,n_l$, $n_1+\ldots+n_l = n$, and any $R_1,\ldots, R_l>0$).


\hfill


For the proof see Section~\ref{_Packing_by_arb_shapes_main_Section_}.


\hfill


\question
\\
Does \ref{_full_packing_of_tori_by_polydisks_Corollary_} hold
for all K\"ahler forms on $T^{2n}$?


\hfill


\remark \label{_Two-dim-case_Remark_}

In dimension 2, volume (that is, symplectic
area) is the only constraint for
symplectic packing of any closed 2-dimen\-si\-onal symplectic
manifold by any shapes. This easily follows from Moser's theorem
\cite{_Moser_} stating that two symplectic forms on a closed surface
are symplectomorphic if and only if their integrals over the surface
are equal. (This explains why we left out the 2-dimensional case in
\ref{_Packing_by_arb_shapes_main_Theorem_} and
\ref{_full_packing_of_tori_by_polydisks_Corollary_}).


\hfill


By the abovementioned Moser's theorem, the closed disk
$B^2(1/\sqrt{\pi})\subset (\R^2, \Omega_2)$ and the square
$[0,1]^2\subset (\R^2, \Omega_2)$ (that both have area 1) have
arbitrarily small symplectomorphic open neighborhoods and therefore
so do the $2n$-dimensional polydisk ${\mathcal
Poly}:=B^2(1/\sqrt{\pi})\times\ldots\times B^2(1/\sqrt{\pi})\subset
(\R^{2n}, dp\wedge dq)$ and the $2n$-dimensional cube ${\mathcal
Cube}:=[0,1]^2\subset (\R^{2n}, dp\wedge dq)$. Therefore $$\nu
(M,\omega, {\mathcal Poly}) = \nu (M,\omega, {\mathcal Cube})$$ for
all $(M,\omega)$. Thus,
\ref{_full_packing_of_tori_by_polydisks_Corollary_} yields the
following corollary.


\hfill


\corollary
\label{_full_packing_of_tori_by_cubes_Corollary_}
\\
Let $\omega$ be a K\"ahler form on $T^{2n}$, $n\geq 2$, and assume
that the cohomology class $[\omega]$ is not proportional to a rational one. Then for any $k\in\N$ the symplectic
manifold $(T^{2n},\omega)$ can be fully packed by $k$ equal cubes.
\endproof


\hfill


\section{Campana simple K\"ahler manifolds}
\label{_Campana_intro_Section_}

\definition

A complex structure on a (closed, connected) manifold $M$, $\dim_\C M>1$,
 is called {\bf Campana simple}, if the union ${\goth U}$
of  all complex subvarieties $Z\subset M$ satisfying $0< \dim_\C Z< \dim_\C M$
has measure\footnote{The measure is defined by means of a volume form on $M$. One can easily see
that if a set is of measure zero with respect to some volume form, then it is of measure zero with respect
to any other volume form.} zero.

The points of $M\backslash {\goth U}$ are called {\bf  Campana-generic}.

If a complex structure $J$ on $M$ is Campana simple, we will also say that
the complex manifold $(M,J)$ is {\bf Campana simple}.


\hfill


\remark

If $J$ is a complex structure of K\"ahler type on a closed connected
manifold $M$, then the union ${\goth U}$ of  all complex
subvarieties $Z\subset M$ satisfying $0< \dim_\C Z< \dim_\C M$
either has measure zero or is the whole $M$.

Indeed, the proper
positive-dimensional complex subvarieties of $(M,J)$ are
parameterized by points of the so-called Douady space $\cD$ of $(M,J)$
which is itself a complex variety, and moreover, there exists a complex variety
$\cX$ and proper holomorphic projections $\pi_1: \cX\to \cD$ and $\pi_2:\cX\to M$
so that for each $a\in\cD$ the set $\pi_2 (\pi_1^{-1} (a))$ is exactly the complex subvariety
of $M$ parameterized by $a$ \cite{_Douady_1966}. By a theorem of Fujiki \cite{_Fujiki:compactness_},
each
irreducible component $\cD'$ of $\cD$ is compact. Thus $\pi_1^{-1} (\cD')\subset \cX$
is a compact complex variety. By Remmert's proper holomorphic mapping theorem \cite{_Remmert_MathAnn_1956_}, since $\pi_2$ is holomorphic and proper, $\pi_2 (\pi_1^{-1} (\cD'))$
is a complex subvariety of $(M,J)$.
Note that  $\pi_2 (\pi_1^{-1} (\cD'))$ is the union of all
the complex subvarieties of $(M,J)$ parameterized by the
points of $\cD'$.
By another theorem of Fujiki \cite{_Fujiki:Douady_}, the space
$\cD$ is second-countable and therefore has at most countably many irreducible components.
Thus,
${\goth U}$ is a union of at most countably many complex subvarieties
of $(M,J)$. If all of these subvarieties are proper, ${\goth U}$ has measure 0, otherwise ${\goth U}=M$.


\hfill


\remark

Campana simple manifolds are non-algebraic.
Indeed, a manifold which admits a globally defined meromorphic function
$f$ is a union of zero divisors of the functions $f-a$, for all
$a\in \C$, and the zero divisor of
$f^{-1}$. Hence, unlike algebraic manifolds, Campana simple manifolds
admit no globally defined meromorphic functions.


\hfill


The following conjecture is due to F.Campana.


\hfill

\conjecture (\cite[Question 1.4]{_Campana:isotrivial_},
\cite[Conjecture 1.1]{_CDV:threefolds_})\\
Let $(M,J)$ be a Campana simple K\"ahler manifold.
Then $(M,J)$ is bimeromorphic to a
hyperk\"ahler orbifold or a finite quotient of a torus.


\hfill


\subsection{Deformations of complex structures}

Let $M$, $\dim_\R M=2n$, be a closed connected manifold and let $J$ be a complex structure on $M$.

Assume $\cX$, $\cB$ are connected -- not necessarily Hausdorff! --
complex manifolds, $t_0\in \cB$ is a marked point and $\cX\to \cB$
is a proper holomorphic submersion whose fiber over $t_0\in \cB$ is
$M$. By Ehresmann's lemma, the fibration $\cX\to \cB$ admits {\it
smooth} local trivializations. The fibers of $\cX\to \cB$ are closed
complex submanifolds of $\cX$ that are diffeomorphic to $M$. Denote
by $M_t$ the fiber over $t\in \cB$  and denote by $J_t$ the complex
structure on $M_t$ induced by the complex structure on $\cX$. Assume
that $J_{t_0} = J$. In such a case we say that $\cX\to (\cB,t_0)$ is
{\bf a smooth deformation of $(M,J)$}.

A smooth local trivialization of $\cX\to \cB$ over a small
neighborhood of $t_0$ in $\cB$, identified with a ball
$B^{2m}\subset \C^m$ (so that $t_0$ is identified with $0\in
B^{2m}$), allows to view all $J_t$, $t\in B^{2m}$, as complex
structures on $M$. Such a family $\{ J_t\}$, $t\in B^{2m}$, $J_0=J$,
will be called {\bf a smooth local deformation of $J$}. (Note that
we work only with deformations with a smooth base).

We say that $J$ {\bf can be approximated by Campana-simple complex
structures in a smooth deformation} if there exists a smooth local
deformation $\{ J_t\}$, $t\in B^{2m}$, $J_0=J$, of $J$ and a
sequence $\{ t_i\}\to 0$ in $B^{2m}$ so that each $J_{t_i}$ is
Campana simple.


\hfill


\subsection{Proof of \ref{_FSP_main_Theorem_}}
\label{pf-of-fsp-main-thm_Section}

The proof of \ref{_FSP_main_Theorem_} is based on the following two claims.


\hfill


\theorem\label{_Campana_complex_structure_existence_Theorem_}\\
(A) Any complex structure of K\"ahler type on $T^{2n}$ can be approximated in a smooth deformation by
complex structures $J$
such that $(T^{2n},J)$ does not admit any proper complex subvarieties
of positive dimension and, in particular, is Campana simple.

\noindent (B) Let $(M,\goh)$, $\dim_\R M\geq 4$, be a closed connected IHS
hyperk\"ahler manifold and let $I$ be a complex structure appearing in $\goh$.
Then $I$ can be approximated in a smooth deformation by
Campana simple complex structures.


\hfill


For the proofs of (A) and (B) see
Sections~\ref{_Campana_simple_complex_structures_torus_Section_},
\ref{_Campana_simple_str_hyperkahler_manifolds_Section_}.


\hfill


\theorem\label{_Campana_full_packing_Theorem_}\\
Let $(M,I,\omega)$ be a K\"ahler manifold and assume that $I$ can be
approximated in a smooth deformation by Campana simple complex structures.
Then the symplectic packings of $(M,\omega)$ by balls are unobstructed.


\hfill


Below we give a sketch of the proof -- for the complete proof (in
fact, of a somewhat stronger result)
see
Section~\ref{_sympl_packing_Campana_simple_Kahler_mfds_Section_}
and, in particular,
\ref{_Campana_simple_FSP_precise_version_Theorem_}.


\hfill


\noindent {\bf Proof of \ref{_FSP_main_Theorem_}.}

In the case
$M=T^2$ the theorem is obvious -- see \ref{_Two-dim-case_Remark_}. In all the other cases the proof follows right
away from \ref{_Campana_complex_structure_existence_Theorem_} and
\ref{_Campana_full_packing_Theorem_}. \endproof



\hfill


\subsection{A sketch of the proof of \ref{_Campana_full_packing_Theorem_}}
\label{_Sketch_of_the_proof_of_the_main_thm_on_full_packing_Section_}

The basic notions used in this sketch will be recalled in further sections.

Let $(M,\omega)$ be as in \ref{_Campana_full_packing_Theorem_}.
Assume we want to show that the symplectic packings of $(M,\omega)$
by $k$ balls are unobstructed. Let $\tM$ be a complex blow-up of $M$ at $k$
Campana-generic points $x_1,\ldots,x_k$. More precisely, we think of
$\tM$ as a fixed smooth manifold -- the connected sum of $M$ with
$k$ copies of $\overline{\C P^n}$ -- so that each complex structure
$I$ on $M$ induces a complex structure $\tI$ on $\tM$ and a
projection $\Pi_I : \tM\to M$ whose fibers over $x_1,\ldots, x_k$
are the exceptional divisors $E_1 (I),\ldots , E_k (I)$ that are
complex submanifolds of $(\tM, \tI)$. The exceptional divisors vary
as the complex structure $I$ varies but their homology classes
remain the same. Denote by $[E_1],\ldots, [E_k] \in H^2 (\tM,\Z)$
the corresponding Poincar\'e-dual cohomology classes. Similarly, the
projection $\Pi_I$ varies with $I$ but the induced map $H^*
(M,\C)\to H^* (\tM,\C)$, that will be denoted by $\Pi^*$, remains
the same.

We are going to use the following theorem of McDuff and Polterovich.


\hfill


\theorem (McDuff-Polterovich, \cite[Prop. 2.1.B and Cor.
2.1.D]{_McD-Polt_}) \label{_McD_P_main_Theorem_sketch_}

Let $M$,
$\dim_\R M = 2n$, be a closed connected manifold equipped with a
K\"ahler form $\omega$. Let $r_1,\ldots, r_k$ be a collection of
positive numbers. Assume there exists a complex structure $I$ of
K\"ahler type on $M$ tamed by $\omega$ and a symplectic form
$\tomega$ on $\tM$ taming $\tI$ so that $[\tomega] = \Pi^* [\omega]
- \pi \sum_{i=1}^k r_i^2 [E_i]$. Then $(M,\omega)$ admits a
symplectic embedding of $\bigsqcup\limits_{i=1}^k B^{2n} (r_i)$.

For sufficiently small $r_1,\ldots, r_k>0$ the
cohomology class $[\tomega] = \Pi^* [\omega] - \pi \sum_{i=1}^k
r_i^2 [E_i]$ is K\"ahler.
\endproof


\hfill


We are also going to use the Demailly-Paun
theorem \cite{_Demailly_Paun_} that says the following: Let $N$ be a closed connected K\"ahler manifold
and let $\hat K(N)$ be the subset of $H^{1,1}(N,\R)$
consisting of all classes $\eta$ such that $\int_Z \eta^{\dim Z}>0$
for any closed complex subvariety $Z\subset N$.
The Demailly-Paun theorem says that the K\"ahler
cone of $N$ is one of the connected components of
$\hat K(N)$.

At the first stage of the proof assume that $(M,I,\omega)$ is a
K\"ahler manifold and $I$ is Campana simple. There are three kinds
of irreducible complex subvarieties of $\tM$:

\begin{description}
\item[(a)] The preimages under $\Pi_I$ of proper complex subvarieties
$Z\subset M$ not containing the points $x_i$, $i=1,\ldots,k$.
\item[(b)] The subvarieties of $E_i (I)$, $i=1,\ldots,k$.
\item[(c)] The manifold $\tM$ itself.
\end{description}

Any real $(1,1)$-cohomology class $\eta$ on $\tM$ can be written as
$$\eta=\Pi^* [\omega] - \pi\sum_{i=1}^k r_i^2 [E_i]$$ for some
$r_1,\ldots,r_k>0$. To see when $\eta$ lies in the Demailly-Paun set
$\hat K(\tM)$ we consider separately the three types of subvarieties
defined above
and using the
Demailly-Paun theorem deduce (see Section~\ref{_sympl_packing_Campana_simple_Kahler_mfds_Section_}
for the precise argument) that the K\"ahler cone of $\tM$ is the
set of all $\eta:=\Pi^* [\omega] - \pi\sum_{i=1}^k r_i^2 [E_i]$ that
satisfy the following conditions:

(A) $r_i >0$ for all $i=1,\ldots,k$,

and

(B) $\int_M \eta^n =\int_M [\omega]^n- \pi^n\sum_{i=1}^k r_i^{2n}>0$.

Combining \ref{_McD_P_main_Theorem_sketch_} with this description of
the K\"ahler cone of $\tM$ we immediately obtain that the symplectic packings of $(M,\omega)$
by balls are unobstructed.

Till this point the argument has been similar to the one in \cite{_LMcDS_}.
Let us now consider the general case when $I$ is only a limit of
Campana simple complex structures in a smooth deformation -- this part of the argument
is new, compared to \cite{_LMcDS_}, and yields a stronger result already when $\dim M =4$:
for instance, it implies that for a product symplectic form $\omega=\omega'\oplus\omega'$ on $T^4=T^2\times T^2$ the symplectic packings of $(T^4,\omega)$ by {\it any} number of balls are unobstructed (in \cite{_LMcDS_}
this was proved only for the packings by one ball).

Namely, assume we want to embed
$\bigsqcup\limits_{i=1}^k B^{2n} (r_i)$, $\Vol
(\bigsqcup\limits_{i=1}^k B^{2n} (r_i)) < \Vol (M,\omega)$,
symplectically into $(M,\omega)$. Let $J$ be a Campana simple
complex structure of K\"ahler type close to $I$ in a smooth deformation (in particular, we
may assume that $\omega$ tames $J$) and let $\tJ$ be the
corresponding complex structure $\tJ$ on $\tM$. We use the
Kodaira-Spencer stability theorem \cite{_Kod-Spen-AnnMath-1960_} to
show that the (1,1)-part $[\omega]_J^{1,1}$ of the cohomology class
$[\omega]$ with respect to $J$ can be represented by a K\"ahler form
$\omega'$ compatible with $J$ and close to $\omega$ so that $\Vol
(\bigsqcup\limits_{i=1}^k B^{2n} (r_i)) < \Vol (M,\omega')$
(see Section~\ref{_Kahler_cone_Subsection} for a rigorous proof). Since the symplectic packings by balls in the Campana simple case are unobstructed, together with
\ref{_McD_P_main_Theorem_sketch_}, this yields that the class
$\alpha:=\Pi^* [\omega]_J^{1,1} - \pi\sum_{i=1}^k r_i^2 [E_i]\in H^2
(\tM,\R)$ is K\"ahler. The next step is crucial: we observe that the
cohomology class $\eta = \Pi^* [\omega] - \pi \sum_{i=1}^k r_i^2
[E_i]$ can be written as $\eta=\alpha + \beta$, where $\beta\in H^2
(\tM,\R)$ is a (2,0)+(0,2) class {\it with respect to $\tJ$} (see Section~\ref{_sympl_packing_Campana_simple_Kahler_mfds_Section_}
for the precise argument). This
implies that $\eta$ is a symplectic class. Moreover, it is not hard
to show (see Section~\ref{_sympl_packing_Campana_simple_Kahler_mfds_Section_} for details) that $\eta$ can be represented by a symplectic form taming
$\tJ$. Applying \ref{_McD_P_main_Theorem_sketch_} we obtain that
$(M,\omega)$ admits a symplectic embedding of
$\bigsqcup\limits_{i=1}^k B^{2n} (r_i)$, as required.


\section{Background from complex geometry}

\subsection{Hodge decomposition}
\label{_Hodge_structures_Section_}

We recall a few basic facts about Hodge structures -- for
details see e.g. \cite{_Voisin-Hodge_}.

Let $M$, $\dim_\R M=2n$, be a closed connected manifold admitting K\"ahler structures.

An almost complex structure $J$ on $M$ acts on the space of
$\R$-linear $\C$-valued differential forms on $M$ (the action is
induced by the pointwise action of $J$ on the tangent spaces of $M$).
This action induces a well-known $(p,q)$-decomposition of the space of
such forms. The $(p,q)$-component of a differential form $\omega$
with respect to this decomposition will be denoted by
$\omega^{p,q}_J$; sometimes, in order to emphasize the dependence of
the decomposition on $J$, we will also say that a form is a {\bf
$(p,q)$-form with respect to $J$}. The action of $J$ on a
$(p,q)$-form is simply the multiplication of the form by $i^{p-q}$.
In particular, all the $(p,p)$-forms are preserved by the action.

For the remainder of this section let $J$ be a complex structure of
K\"ah\-ler type on $M$. Then, by the famous theorem of Hodge, saying
that any $(p,q)$-form is cohomologous to a unique harmonic
$(p,q)$-form, the $(p,q)$-de\-com\-po\-si\-tion on the space of
differential forms induces a $(p,q)$-decomposition on $H^* (M,\C)$,
called {\bf the Hodge decomposition}. Although the proof of the
existence of the Hodge decomposition involves the whole K\"ahler
structure of which $J$ is a part, one can show (see e.g. \cite[Vol.
1, Prop. 6.11]{_Voisin-Hodge_}) that, in fact, the Hodge
decomposition depends only on the complex structure $J$. We
will denote it by
$$H^* (M,\C) = \bigoplus\limits_{p,q} H^{p,q}_J (M,\C).$$
We will also say that a (non-zero) cohomology class in $H^{p,q}_J (M,\C)$ is
a {\bf $(p,q)$-class (with respect to $J$)}.

Set
$$H^{p,q}_J (M,\Q):= H^{p,q}_J (M,\C)\cap H^{p+q} (M,\Q),$$
$$H^{p,q}_J (M,\R):= H^{p,q}_J (M,\C)\cap H^{p+q} (M,\R).$$

The following proposition is also well-known -- see e.g. \cite[Vol. 1, Sec. 11]{_Voisin-Hodge_}.


\hfill


\proposition
\label{_fund_class_cplx_subvariety_Proposition_}\\
Let $J$ be a complex structure of K\"ahler type on a closed connected manifold $M$, $\dim_\R M = 2n$.
If $L\subset (M,J)$, $\dim_\R L = 2l$,  is a closed complex subvariety, then $L$ represents a well-defined fundamental integral
homology class of degree $2l$; denote its Poincar\'e-dual cohomology class by $[L]\in H^{2n-2l} (M,\Q)$. Then
$[L]\in H^{n-l,n-l}_J (M,\Q)$ and $[L]\neq 0$.
\endproof


\hfill


Let $\cX\to (\cB,t_0)$ be a smooth deformation of $(M,J)$ over a
connected base $\cB$ which is trivial as a smooth fibration. Assume
that for each $t\in\cB$ the complex structure $J_t$ on the fiber
$M_t$ of $\cX\to (\cB,t_0)$ is of K\"ahler type. Note that there is a canonical (that is,
independent of trivializations) identification of the
homology/cohomology of each fiber $M_t$ with the homology/cohomology
of $M$. Thus for each $t\in\cB$ we have a well-defined Hodge
decomposition of $H^* (M,\C)$ defined by $J_t$.


\hfill


\definition
\\
Let $a\in H^{2p} (M,\Q)$, $a\neq 0$. The {\bf Hodge locus of $a$ for
the smooth local deformation $\cX\to (\cB,t_0)$} is the set of $t\in
\cB$ such that $a\in H_{J_t}^{p,p} (M,\Q)$.


\hfill


The following proposition can be found e.g. in \cite[Vol. 2, Lem. 5.13]{_Voisin-Hodge_}.


\hfill


\proposition
\label{_set_cplx_str_for_which_given_class_p-p_complex_subvariety_Proposition_}\\
The Hodge locus of $a\in H^{2p} (M,\Q)$, $a\neq 0$, is a complex subvariety of $\cB$. \endproof


\hfill


\subsection{Complex structures tamed by symplectic forms}
\label{_complex_structures_tamed_by_sympl_forms_Section_}

\definition
\\
We say that a differential 2-form $\omega$ {\bf tames} an almost complex
structure $J$ if $\omega (v,Jv) >0$ for any non-zero tangent vector $v$.


\hfill


Clearly, a closed differential 2-form $\omega$ taming an almost complex structure $J$ is automatically symplectic.


\hfill

\proposition\label{_tamed+2,0_Proposition_}
\\
(A) A differential 2-form $\omega$ tames an almost complex structure
$J$ if and only if so does $\omega_J^{1,1}$.

\noindent
(B) A differential 2-form $\omega$ is compatible with an almost complex
structure $J$ if and only if $\omega$ tames $J$ and $\omega =
\omega_J^{1,1}$.

\noindent
(C) Assume $\omega$ is a symplectic form on a manifold $M$ taming
an almost complex structure $J$. Let $\eta$ be a closed
(real-valued) $(2,0)+(0,2)$-form with respect to $J$. Then $\omega+\eta$ is also a
symplectic form.


\hfill


\noindent {\bf Proof:}

The $(2,0)+(0,2)$-forms (with respect to $J$) are
exactly the 2-forms anti-invariant under the action of $J$,
while the $(1,1)$-forms are exactly the invariant ones. In
particular, if $\eta$ is a $(2,0)+(0,2)$-form, then $\eta (v,Jv)=0$
for all $v$. The claims (A), (B) and (C) now follow easily.
\endproof


\subsection{Symplectic and K\"ahler cones}
\label{_Kahler_cone_Subsection}


Let $M$ be a closed connected smooth manifold.

A cohomology class $a\in H^2 (M,\R)$ is called {\bf symplectic} if
it can be represented by a symplectic form. The set of all symplectic classes $a\in H^2 (M,\R)$ is called {\bf the
symplectic cone} of $M$.

Assume $J$ is a complex structure on $M$.
A cohomology class $a\in H^2 (M,\R)$ is called
{\bf K\"ahler}, or {\bf K\"ahler with respect to $J$}, if it can
be represented by a K\"ahler form compatible with $J$. The set of all K\"ahler classes in $H^2 (M,\R)$
is a convex cone $\Kah (M,J)$, called the {\bf K\"ahler cone} of
$(M,J)$. Clearly, it is a subset of the symplectic cone of $M$.


\hfill

\theorem(A version of Kodaira-Spencer
stability theorem)
\label{_Kodaira_stabi_Theorem_}\\
Let $(M,J,\omega)$ be a closed K\"ahler manifold,
and let $\{ J_t\}$, $t\in
B$, $J_0=J$, be a smooth local deformation of $J$.
Then there exists a neighborhood of $U\subset B$ of zero in $B$
such that the complex structure $J_t$ on $M$ is of K\"ahler type and
$[\omega]^{1,1}_{J_t}\in \Kah (M,J_t)$ for all $t\in U$.
Moreover, the class $[\omega]^{1,1}_{J_t}$, $t\in U$, depends smoothly on $t$.
\endproof


\hfill


\noindent {\bf Proof:}

The Kodaira-Spencer stability theorem \cite{_Kod-Spen-AnnMath-1960_} states that there exists an open neighborhood
$U\subset B$ of zero such that $J_t$ is of K\"ahler type for any $t\in U$ and $\omega$
can be
extended to a smooth family $\{ \omega_t\}$, $t\in U$, $\omega_0=\omega$, of K\"ahler
forms on the complex manifolds $(M, J_t)$.

For $t\in U$ denote by $\Omega_t$ the $(1,1)$-component of the
unique harmonic 2-form on the closed K\"ahler manifold $(M,J_t,\omega_t)$
representing the cohomology class $[\omega]\in H^2 (M,\R)\subset H^2
(M,\C)$. Then $\Omega_t$ is a real closed 2-form on $M$ of type
$(1,1)$ (with respect to $J_t$) and $[\Omega_t] =
[\omega]^{1,1}_{J_t}$ for all $t\in U$.

Note that $\Omega_0$ and $\omega$ are cohomologous closed
real-valued $(1,1)$-forms on $(M,J)$. Therefore, by the
$\partial\bar{\partial}$-lemma (see e.g. \cite[Vol.1, Prop. 6.17]{_Voisin-Hodge_}), $\omega = \Omega_0 + \partial_J
\bar{\partial}_J F$ for a smooth function $F: M\to \R$.

Consider now the forms $\alpha_t:= \Omega_t + \partial_{J_t}
\bar{\partial}_{J_t} F$, $t\in U$. The form $\alpha_t$ is a closed
$(1,1)$-form with respect to $J_t$ which is cohomologous to
$\Omega_t$. Thus $[\alpha_t] = [\Omega_t] = [\omega]^{1,1}_{J_t}$.

Since $\Omega_t$ depends
smoothly on $t$, so does $\alpha_t$ and therefore so does its cohomology class $[\alpha_t]$.

Since the condition on a complex structure to be tamed by a
symplectic form is open and since $\omega=\alpha_0$ tames $J=J_0$,
we can assume without loss of generality that $U$ is sufficiently
small so that the form $\alpha_t$ tames $J_t$ for all $t\in U$. By
\ref{_tamed+2,0_Proposition_}, part B, this means that $\alpha_t$ is
K\"ahler on $(M,J_t)$. Hence, $[\alpha_t] = [\omega]^{1,1}_{J_t}$ lies in
$\Kah(M,J_t)$ and depends smoothly on $t$.
\endproof


\hfill


\section{Campana simple complex structures on a torus}
\label{_Campana_simple_complex_structures_torus_Section_}

Let $M=T^{2n}$, $n\geq 2$.
As above, let $\pi: \R^{2n}\to \R^{2n}/\Z^{2n}=T^{2n}$ be the standard projection.


\hfill


\proposition
\label{_Any_complex_structure_of_Kahler_type_on_torus_linear_Proposition_}\\
Any complex structure of K\"ahler type on $T^{2n}$ is biholomorphic
to a linear complex structure. Any K\"ahler form on $T^{2n}$ is
symplectomorphic to a linear symplectic form.
The biholomorphism and the symplectomorphism can chosen so that they act on $H^* (T^{2n})$ by identity.


\hfill


\noindent {\bf Proof.}

Let $(J,\omega)$ be a K\"ahler structure on $T^{2n}$. The
Albanese
map (see e.g. \cite[Vol. 1, Def. 12.10 and Thm.
12.15]{_Voisin-Hodge_} defines a biholomorphism $f: (T^{2n},J)\to
(\C^n/\Gamma, I)$, where $I$ is the standard complex structure on
$\C^n/\Gamma$ and $\Gamma\subset \C^n$ is a lattice.

Denote by $\pi_\Gamma: \C^n\to \C^n/\Gamma$ the natural projection
and let $F:\R^{2n}\to\C^n$ be an $\R$-linear isomorphism of
vector spaces mapping $\Z^{2n}$ to $\Gamma$. Then $F$ covers a
diffeomorphism $\bar{f}: T^{2n} = \R^{2n}/\Z^{2n}\to \C^n/\Gamma$,
that is, $\bar{f}\circ \pi = \pi_\Gamma\circ F$.
One easily sees that $F$ can be chosen in such a way that $\bar{f}$ induces the same map $H^* (\R^{2n}/\Z^{2n})\to H^* (\C^n/\Gamma)$
as $f$.

Clearly, $\bar{f}^* I$ is a linear complex structure on
$T^{2n}=\R^{2n}/\Z^{2n}$. Therefore $\bar{f}^{-1}\circ f$ is a
diffeomorphism mapping $J$ to a linear complex structure on
$T^{2n}$
and acting as identity on $H^* (T^{2n})$. This proves the first claim of the proposition.

Let us prove the second claim. The symplectic form
$\zeta:=(f^{-1})^* \omega$ on the torus $\C^n/\Gamma$ is compatible
with $I$ and is cohomologous to a linear 2-form $\eta$ (i.e.
$\pi_\Gamma^*\eta$ is a linear symplectic form on $\C^n$). By
\ref{_tamed+2,0_Proposition_}, $\zeta = \zeta_I^{1,1}$, i.e. $\zeta$
is a $(1,1)$-form. Since $\eta$ is linear, it is harmonic, and since
it is cohomologous to $\zeta$, it is also a $(1,1)$-form. Therefore,
by \ref{_tamed+2,0_Proposition_}, $\eta$ is symplectic.

Now note that $\zeta$ and $\eta$ are two cohomologous symplectic
forms on $\C^n/\Gamma$ that are compatible with the same complex structure $I$. Therefore
$\zeta$ and $\eta$ can be connected by a straight path $t\zeta +
(1-t)\eta$, $t\in [0,1]$, of closed cohomologous $2$-forms; all
these forms tame $I$ and hence are symplectic. By Moser's theorem
\cite{_Moser_}, it implies that the symplectic forms $\zeta$ and
$\eta$ on $\C^n/\Gamma$ are
mapped into each other by a diffeomorphism of $\C^n/\Gamma$ isotopic to the identity.
Hence the symplectic
forms $\omega = f^* \zeta$ and $\bar{f}^* \eta$ on
$T^{2n}=\R^{2n}/\Z^{2n}$ are also
mapped into each other by a diffeomorphism of $T^{2n}$ isotopic to the identity.
Clearly, $\bar{f}^*
\eta$ is a linear symplectic form because $\pi^* \bar{f}^* \eta =
(\pi_\Gamma\circ F)^* \eta$. Thus, $\omega$ is
mapped into a
linear symplectic form on $T^{2n}=\R^{2n}/\Z^{2n}$ by a diffeomorphism acting on $H^* (T^{2n})$ by identity.
This finishes the
proof of the proposition.
\endproof


\hfill


\proposition
\label{_moduli_space_of_lin_cplx_str_torus_Proposition_}\\
The space of linear orientation-preserving complex structures on $\T^{2n}$ (which we identify with the space of orientation-preserving linear
complex structures on $\R^{2n}$) is a connected manifold $\Comp_l$ that can be equipped with a complex structure.


\hfill


\noindent {\bf Proof.}

It is well-known (see e.g. \cite{_McD-Sal-intro_}, Prop. 2.48) that the space of the orienta\-tion-preserving linear
complex structures on $\R^{2n}$ can be identified with
$$\Comp_l:=GL^+ (2n,\R)/GL(n,\C)$$ which is a connected smooth manifold.
Here the group $GL (n,\C)$ is embedded in $GL (2n,\R)$ by the map
\[
C\mapsto \left(
           \begin{array}{cc}
             \textrm{Re}\ C & -\textrm{Im}\ C \\
             \textrm{Im}\ C & \textrm{Re}\ C \\
           \end{array}
         \right)
\]
The tangent space $T_J \Comp_l$ to $\Comp_l$ at a point $J\in
\Comp_l$ is given by all $R\in GL^+ (2n,\R)$ satisfying $RJ+JR=0$.
The complex structure on $T_J \Comp_l$ is given by $R\mapsto JR$.
Since $\Comp_l$ is identified with the symmetric space $GL^+
(2n,\R)/GL(n,\C)$, the
 almost complex structure on $\Comp_l$ is integrable
\cite{_Besse:Einst_Manifo_}.
\endproof


\hfill


Consider the projection $T^{2n}\times \Comp_l\to \Comp_l$. Denote by
$I$ the complex structure on $\Comp_l$. Equip $T^{2n}\times \Comp_l$
with an almost complex structure which is defined at a point
$(x,J)\in T^{2n}\times \Comp_l$ as $J\oplus I$ with respect to the
obvious splitting of $T_{(x,J)} (T^{2n}\times\Comp_l)$. One easily
checks that this almost complex structure is integrable and thus
$T^{2n}\times \Comp_l\to \Comp_l$ is a smooth deformation of
$(T^{2n}, J)$ for any $J\in \Comp_l$.

Given $A\in H^{2p} (T^{2n},\Q)$, $A\neq 0$, $0<p<n$, define
$\cT_A\subset \Comp_l$ as the Hodge locus of $A$ with respect to the
smooth deformation $T^{2n}\times \Comp_l\to \Comp_l$.


\hfill

\proposition
\label{_cT_A_is_proper_complex_subvariety_Proposition_}\\
The set $\cT_A$ is either empty or a proper complex subvariety of
$(\Comp_l, I)$.


\hfill


\noindent {\bf Proof.}

Given a complex structure $J$ of K\"ahler type on $T^{2n}$, denote by
$\bJ$ the linear
complex structure on $\R^{2n}$ which is the lift of $J$ under $\pi: \R^{2n}\to \R^{2n}/\Z^{2n}=T^{2n}$
(see \ref{_Any_complex_structure_of_Kahler_type_on_torus_linear_Proposition_}).

For any cohomology class $A\in H^* (T^{2n}, \R)$ there is a unique
closed linear form $\omega_A$ on $T^{2n}$ that represents $A$. Let $\bomega_A := \pi^* \omega_A$
be the corresponding linear form on $\R^{2n}$.

Now let $A\in H^{2p} (T^{2n},\R)$, $A\neq 0$, $0<p<n$.
It follows from \ref{_set_cplx_str_for_which_given_class_p-p_complex_subvariety_Proposition_}
that $\cT_A$ is a complex subvariety of $\Comp_l$. Let us prove that the complex subvariety $\cT_A$ is proper.

Indeed, assume by contradiction that it is not proper. Then, since $\Comp_l$
is connected, $\cT_A$ has to coincide with $\Comp_l$ -- in other words, $A\in H^{p,p}_J (T^{2n},\R)$ for any $J\in \Comp_l$. Hence,
the linear form $\bomega_A$
is $\bJ$-invariant for any orientation-preserving linear complex structure $\bJ$ on $\R^{2n}$ and thus
it is preserved by the subgroup $G\subset SL (2n,\R)$ generated by such $\bJ$ (each such $\bJ$ lies in $SL (2n,\R)$). Since the set
of orientation-preserving linear complex structures is conjugacy-invariant in $SL (2n,\R)$, $G$ is a normal subgroup of $SL (2n,\R)$.
But $SL (2n,\R)$ is a simple Lie group and therefore its normal subgroup is either contained in its center (equal to $\{ Id, -Id\}$) or
coincides with the whole $SL (2n,\R)$ (see e.g. \cite{_Ragozin_PAMS72_}). Clearly, the first option does not hold for $G$ and therefore
$G=SL (2n,\R)$. Therefore the exterior form $\bomega_A$ is $SL (2n,\R)$-invariant. But the only $SL (2n,\R)$-invariant exterior forms
on $\R^{2n}$ are scalar multiples of a volume form\footnote{This can be easily seen by considering the action of diagonal matrices of the form
\[\textrm{diag}\, \{ 1,\ldots, 1, \lambda,1,\ldots, 1, 1/\lambda, 1, \ldots, 1\}\]
lying in $SL (2n,\R)$ on exterior forms.} or forms of degree $0$
(that is, constants), while $\bomega_A$ is of degree $0<2p<2n$. Thus
we have obtained a contradiction. Hence, $\cT_A$ is a proper
subvariety of $\Comp_l$ and the proposition is proved. \endproof


\hfill


The
following theorem is a precise formulation of
\ref{_Campana_complex_structure_existence_Theorem_} in the case of a
torus.


\hfill

\theorem
\label{_Campana_simple_str_dense_torus_Theorem_}\\
Let $\Xi$ be the set of linear complex structures $J\in \Comp_l$
such that $(T^{2n},J)$ admits no proper complex subvarieties of
positive dimension. Then $\Xi$ is dense in $\Comp_l$.


\hfill


In fact, as it follows from the proof below, $\Xi$ is not only dense
in $\Comp_l$ but of ``full measure".
\ref{_cT_A_is_proper_complex_subvariety_Proposition_} and
\ref{_Campana_simple_str_dense_torus_Theorem_}, which is easily
deduced from \ref{_cT_A_is_proper_complex_subvariety_Proposition_},
seem to be well known but we have not been able to find them in the
literature.


\hfill
\eject


\noindent {\bf Proof of
\ref{_Campana_simple_str_dense_torus_Theorem_}.}

Let $\cC\subset\Comp_l$ be the set of linear complex structures $J$
for which $H^{p,p}_J (T^{2n},\Q)=0$ for all $0<p<n$. By
\ref{_fund_class_cplx_subvariety_Proposition_}, each proper
positive-dimensional compact complex subvariety of a K\"ahler
manifold carries a non-zero integral fundamental class whose
Poincar\'e-dual cohomology class is of Hodge type $(p,p)$ for some
$0<p<n$. Therefore $\cC\subset \Xi$ and it is enough to show that
$\cC$ is dense in $\Comp_l$. Note that $\cC$ is the complement of
$\bigcup_A \cT_A$ in $\Comp_l$, where the union is taken over all
$A\in H^{2p} (T^{2n},\Q)$, $A\neq 0$, for all $0<p<n$. But the
latter union is a countable union of proper complex subvarieties and
therefore (by Baire's theorem) its complement in $\Comp_l$ -- that is,
$\cC$ -- is dense. This finishes the proof. \endproof


\hfill


\noindent {\bf Proof of \ref{_Campana_complex_structure_existence_Theorem_} in the torus case.}

By
\ref{_Any_complex_structure_of_Kahler_type_on_torus_linear_Proposition_},
it suffices to prove the claim of
\ref{_Campana_complex_structure_existence_Theorem_} for an arbitrary
linear complex structure $J$ on $T^{2n}$. Since $T^{2n}\times
\Comp_l\to \Comp_l$ yields a smooth deformation of $J$, the claim
follows immediately from
\ref{_Campana_simple_str_dense_torus_Theorem_}.
\endproof

\section{Campana simple complex structures on an IHS hyper\-k\"ahler manifold}
\label{_Campana_simple_str_hyperkahler_manifolds_Section_}

First, let us recall a few relevant results about topology and
deformation theory of hyperk\"ahler manifolds -- for more details
see \cite{_V:Torelli_}, \cite{_Huybrechts:lec_}, cf. \cite{_Douady_,
_Kodaira_, _Catanese:moduli_}.

Let $(M,\goh)$, $\dim_\R M =
4n$, be a closed connected IHS hyperk\"ahler manifold, $\goh=\{ I_1,
I_2,I_3, \omega_1,\omega_2,\omega_3\}$.

\subsection{Bogomolov-Beauville-Fujiki form}
\label{_BBF_form_Section_}

The content of this section will be used only in Section~\ref{_Ergodic_Subsection_}. We put it here as it belongs
to the general theory of hyperk\"ahler manifolds.


\hfill

\theorem (Fujiki formula, \cite{_Fujiki:HK_})
\label{_BBF_Fujiki_Theorem_}\\
Let $M$ be an IHS hyperk\"ahler manifold as above.

Then there exists a
primitive integral quadratic form $q$ on $H^2(M,\R)$,
and a positive rational number $\kappa$ so that
$\int_M \eta^{2n}=\kappa q(\eta,\eta)^n$ for each
$\eta\in H^2(M,\R)$.

Moreover, $q$ has signature
$(3, b_2-3)$, where $b_2 = \dim_\R H^2 (M,\R)$,
where $b_2
>
3$.
\endproof

\hfill

\definition
\\
The constant $\kappa$ is called {\bf Fujiki constant},
and the quadratic form $q$ is called
{\bf Bo\-go\-mo\-lov-Beauville-Fujiki form}. It is defined
by the Fujiki formula uniquely, up to a sign.

To fix the sign, note that the complex-valued form $\Omega =
\omega_2 + \sqrt{-1}\omega_3$ is a closed, non-degenerate,
holomorphic $(2,0)$-form with respect to $I_1$. Let $\bar\Omega$ be
the complex conjugate $2$-form with respect to $I_1$.

The sign of $q$ is determined
from the following formula (Bogomolov, Beauville \cite{_Beauville_}):
\begin{align*}  c q(\eta,\eta) &=
   (n/2)\int_M \eta\wedge\eta  \wedge \Omega^{n-1}
   \wedge \bar{\Omega}^{n-1} -\\
 &-(1-n)\left(\int_M \eta \wedge \Omega^{n-1}\wedge \bar
   \Omega^{n}\right) \left(\int_M \eta \wedge \Omega^{n}\wedge \bar{\Omega}^{n-1}\right),
\end{align*}
where $c>0$ is a positive constant.

\subsection{Teichm\"uller space for hyperk\"ahler manifolds}
\label{_Teichmuller_space_hyperkahler_manifolds_Subsection_}

Denote by $\Comp_h$ the set
of complex structures of hyperk\"ahler type
on $M$. (Recall that we view complex structures as tensors on
$M$, that is, as integrable almost complex structures).

The set $\Comp_h$ can be
equipped with $C^\infty$-topology in a standard way (see e.g. \cite{_Hamilton:Nash_}) so that
a sequence of complex structures converges in the topology
if it converges uniformly with all derivatives.

Let $\Diff_0\subset \Diff^+$ be the group of smooth isotopies of
$M$. Note that $\Diff_0$ acts on $\Comp_h$.


\hfill


\definition
\\
The quotient topological
space $\Teich := \Comp_h/\Diff_0$ is called {\bf the Teichm\"uller
space of $M$}. For $J\in \Comp_h$ we denote by $[J]$ the corresponding point in
$\Teich$.


\hfill


For any $a,b,c\in\R$, $a^2+b^2+c^2=1$, the tensor $aI_1+bI_2+cI_3$
is a complex structure of hyperk\"ahler type on $M$. We will call
any such  $aI_1+bI_2+cI_3$ {\bf an induced complex structure on $M$}
or {\bf a complex structure induced by $\goh$}. Clearly, all induced
complex structures lie in the same connected component of $\Comp_h$
whose image under $\Comp_h\to \Teich = \Comp_h/\Diff_0$ will be
denoted by $\Teich_0$.

By \cite{_Bogomolov:defo_},
$\Teich_0$ admits the structure of a smooth
(possibly, non-Hausdorff) complex
manifold.
More
precisely, a combination of the fundamental theorems of
Kuranishi \cite{_Kuranishi_} and Bogomolov-Tian-Todorov
\cite{_Bogomolov:defo_, _Tian_,_Todorov_} shows that the Kuranishi
space of $(M,J_0)$ (the base of a certain universal smooth local
deformation of a complex structure defined in \cite{_Kuranishi_}) is
a smooth complex manifold as long as the canonical bundle of
$(M,J_0)$ is trivial, which, of course, is true in the hyperk\"ahler
case. Moreover (see e.g. \cite{_Catanese:moduli_}), there exists a
homeomorphism between a neighborhood
$U_{[J_0]}$
of $[J_0]$ in $\Teich_0$
and an open subset of the Kuranishi space of $(M,J_0)$. The homeomorphism maps
each $[J]\in \Teich_0$ to the point
$x$ of the Kuranishi space
corresponding to a complex structure on $M$ representing $[J]$.
The
smooth and the complex structures on $U_{[J_0]}$ are the pullbacks of the
corresponding structures on the Kuranishi space. The pullback of
the universal smooth local deformation over the Kuranishi space
under the homeomorphism, together with a smooth trivialization of the universal fibration over a neighborhood of $x$ in the Kuranishi space, induces a smooth deformation $M\times
U_{[J_0]}\to U_{[J_0]}$.

The smooth trivialization of the universal fibration over a neighborhood of $x$ is defined up to the action of $\Diff_0$ on the fibers, meaning that different smooth trivializations yield complex structures on the fibers of $M\times
U_{[J_0]}\to U_{[J_0]}$ that differ by the action of $\Diff_0$ on $M$.
Since
$\Diff_0$ acts trivially on the homology of $M$, different complex
structures representing the same point in $\Teich_0$ define the same
Hodge structure on $H^* (M,\C)$. Thus, the Hodge decomposition on the
cohomology of a fiber of $M\times U\to U$ does not
depend on any of the choices).

Thus, for each  $A\in H^{2p} (M,\Q)$, $A\neq 0$, we have a well-defined set $\hcT_A\subset \Teich_0$ formed by all $[J]\in\Teich_0$ such that $A\in H^{p,p}_J (M,\Q)$ for some (or, equivalently, any) complex structure $J$ representing $[J]$. By the discussion above and by
\ref{_set_cplx_str_for_which_given_class_p-p_complex_subvariety_Proposition_},
the intersection of $\hcT_A$ with each $U_{[J]}$, $[J]\in\Teich_0$, is a complex subvariety of $U_{[J]}$. Therefore
$\hcT_A$ is a complex
subvariety of $\Teich_0$.

\subsection{Trianalytic subvarieties}
\label{_trianalytic_varieties_Section_}

\definition (\cite{_Verbitsky:Symplectic_II_})
\label{_trianalytic_Definition_}
\\
 A closed subset $X\subset M$ of the
hyperk\"ahler manifold $(M,\goh)$ is called {\bf trianalytic} (with respect to $\goh$), if $X$
is a complex
subvariety of $(M,I)$ for each complex structure $I$ induced by
$\goh$.

\hfill

\theorem
\label{_G_M_invariant_implies_trianalytic_Theorem_}\\
Let $A\in  H^* (M,\Q)$, $A\neq 0$. Then $A=[N]$ for a trianalytic $N\subset M$ if and only if $A$ is invariant with respect to the action of any
induced complex structure on $H^* (M,\Q)$.

\hfill


\noindent {\bf Proof:}

This is Theorem 4.1 of
\cite{_Verbitsky:Symplectic_II_}.
\endproof

\hfill

\definition \label{_generic_manifolds_Definition_}
\\
Let $J$ be a complex structure of K\"ahler type on $M$.
We say that $J$ is {\bf of general type}
with respect to the hyperk\"ahler structure $\goh$ on $M$, if all the elements of the group
\[ \bigoplus\limits_p H_J^{p,p}(M,\Q)\subset H^*(M,\Q)\]
are invariant under the action of any complex structure induced by $\goh$.

\hfill

\ref{_G_M_invariant_implies_trianalytic_Theorem_} has the following
immediate corollary:


\hfill

\corollary \label{_triana_gene_type_Corollary_}\\
Assume that a complex structure $J$ of K\"ahler type on $M$ is
of general type with respect to $\goh$. Let $L\subset (M,J)$ be a closed
complex subvariety. Then $L$ is trianalytic with respect to
$\goh$.
\endproof



\hfill


Deformations of trianalytic subvarieties were
studied in \cite{_Verbitsky:Deforma_} where the following theorem  was
proved.


\hfill

\theorem(\cite[Theorem 9.1, Theorem 1.2]{_Verbitsky:Deforma_})
\label{_defo_prod_Theorem_}\\
Let $Z\subset (M,\goh)$ be a compact trianalytic subvariety. Then all deformations
of $Z$ are isometric\footnote{It means that $Z$ and any subvariety $Z'$ of $M$ obtained by a deformation of $Z$
  in the class of complex subvarieties of $M$ are isomorphic as complex varieties and the isomorphism respects the Riemannian
  metrics  obtained by the restriction of the Riemannian metric induced by $\goh$ to $Z$ and $Z'$.}  to $Z$, and their union is locally isometric
to a product $Z\times W$, where $W$ is a
hyperk\"ahler variety.
\endproof

\hfill


This yields the following corollary,
used further on in this paper.


\hfill

\corollary \label{_campana_simple_hk_Corollary_}\\
Assume that the hyperk\"ahler manifold $(M,\goh)$ is IHS and $J$ is a complex structure on $M$
which is of general type with respect to
$\goh$.
Then the union of all proper positive-dimensional complex subvarieties of $(M,J)$ (by \ref{_triana_gene_type_Corollary_} all such subvarieties are trianalytic)
has measure 0.


\hfill


\noindent {\bf Proof:}

Let ${\goth J}$ be the set of all deformation
classes of proper positive-dimensional subvarieties of $(M,J)$.
For each
$\alpha \in {\goth J}$ the union of all subvarieties belonging to $\alpha$
is either a proper complex analytic subvariety of $M$
\cite{_Fujiki:compactness_}
or the whole $M$. However, by \ref{_defo_prod_Theorem_}, in the latter case the hypek\"ahler manifold
$M$ is (locally) a direct product
of two hyperk\"ahler manifolds and therefore is not IHS \cite{_Besse:Einst_Manifo_}, contrary to our assumptions. Thus the former case holds
and therefore the union of all subvarieties belonging to each $\alpha\in {\goth J}$ is of measure zero.
Since ${\goth J}$ is countable \cite{_Fujiki:Douady_}, the union of all proper positive-dimensional complex subvarieties of $(M,J)$
has measure 0.
\endproof

\subsection{Density of Campana simple complex
structures on a hyperk\"ahler manifold}
\label{_density_of_Campana_simple_complex_structures_hyperkahler_Section_}

Note that the action of $\Diff_0$ on $\Comp_h$ maps Campana simple
complex structures to Campana simple complex structures.


\hfill

\theorem
\label{_Campana_simple_str_dense_hyperkahler_Theorem_}\\
Assume that the hyperk\"ahler manifold $M$ is IHS. Then
the set
$$\Camp:=\{ [J]\in\Teich_0\ |\ J\ \textrm{is Campana simple}\}$$
is dense in $\Teich_0$.
Consequently, the set of Campana simple complex structures of hyperk\"ahler type on $M$ is dense in $\Comp_h$.


\hfill


\noindent {\bf Proof of
\ref{_Campana_simple_str_dense_hyperkahler_Theorem_}.}

Define $\hcC\subset\Teich_0$
as $\hcC:=\bigcup_A \hcT_A$,
where the union is taken over all  non-zero $A\in H^* (M,\Q)$ such that $\hcT_A$ is a proper subvariety of $\Teich_0$.
Then $\hcC$ is a countable union of proper complex subvarieties and therefore (by Baire's theorem)
its complement in $\Teich_0$ is dense. Let us show that
$\Teich_0 \setminus \hcC  \subset \Camp$ -- this would prove the theorem.

Indeed, let $J$ be a complex structure of K\"ahler type on $M$ such
that $[J]\in \Teich_0\setminus\hcC$. Let us prove that $J$ is Campana
simple.

Consider the set of $A\in H^* (M,\Q)$, $A\neq 0$, $0< {\rm deg}\,
A<2n={\dim_\R M}/2$, such that $J\in \hcT_A$.

If this set is empty, then \ref{_fund_class_cplx_subvariety_Proposition_} implies that $J$ is Campana simple.

If the set is not empty, pick an arbitrary $A\in H^{p,p}_J (M,\Q)$,
$A\neq 0$, $0<p<2n$,
such that $J\in \hcT_A$.
Recall that $\hcT_A$ is a complex subvariety of $\Teich_0$. By the choice of
$J$, this complex subvariety is non-proper, which means that $\hcT_A = \Teich_0$, since
$\Teich_0$ is connected. In other words, $A\in H^{p,p}_{J'} (M,\Q)$
for all complex structures $J'$ such that $[J']\in \Teich_0$ and, in
particular, $A\in H^{p,p}_{J'} (M,\Q)$ for all induced complex
structures $J'$. Since it holds for arbitrary $A\in H^{p,p}_J
(M,\Q)$, $A\neq 0$, $0<p<2n$,
such that $J\in \hcT_A$,
we get that $J$ is of general type
with respect to $\goh$. By \ref{_triana_gene_type_Corollary_}, any
complex subvariety of $(M,J)$ is trianalytic. By
\ref{_campana_simple_hk_Corollary_}, the union of all such
subvarieties is of measure zero, which means that $J$ is Campana
simple. This finishes the proof. \endproof


\hfill

\noindent {\bf Proof of \ref{_Campana_complex_structure_existence_Theorem_} in the hyperk\"ahler case.}

As we saw above, $M\times \Teich_0\to\Teich_0$ is a smooth
deformation of any complex structure appearing in $\goh$. Therefore
the claim of \ref{_Campana_complex_structure_existence_Theorem_} in
the hyperk\"ahler case follows immediately from
\ref{_Campana_simple_str_dense_hyperkahler_Theorem_}.
\endproof


\section{Symplectic packing for Campana simple K\"ahler manifolds}
\label{_sympl_packing_Campana_simple_Kahler_mfds_Section_}


\subsection{Blow-ups and McDuff-Polterovich theorem}
\label{_blow_ups_McDuff_Polterovich_Section_}


The blow-up operation can be performed both in complex and symplectic categories.

Since we are going to compare blow-ups of the
same smooth manifold with different complex structures, let us take the following point of view on the (simultaneous)
complex blow-ups of a complex manifold $M$ at $k$ points. From this point on we fix $k\in\N$.

We first define a smooth manifold $\tM$ as a connected sum of $M^{2n}$ with $k$ copies of $\overline{\C P^n}$.
Any complex structure $I$ on $M$ defines uniquely, up to a smooth isotopy,
a complex structure $\tI$ on $\tM$ and complex submanifolds $E_1 (I), \ldots, E_k (I)\subset (\tM, \tI)$
so that there exists a diffeomorphism between $\tM$ and the complex blow-up
$\tM'$ of $(M,I)$ at $k$ points identifying $\tI$ with the
canonical complex structure on $\tM'$ and $E_1 (I), \ldots, E_k (I)$ with the $k$ exceptional divisors in $\tM'$.
We will call $E_1 (I), \ldots, E_k (I)$ {\bf the exceptional divisors defined by $I$}.
The canonical projection $\tM'\to M$ is then identified with a projection $\Pi_I: \tM\to M$.
If $J$ is another complex structure on $M$, we get another complex structure $\tJ$ on $\tM$
with another projection $\Pi_J: \tM\to M$ which is smoothly isotopic to $\Pi_I$ and
therefore induces the same map on cohomology which is independent of the
complex structure and will be denoted by $\Pi^*$. The exceptional divisors
$E_1 (J), \ldots, E_k (J)$ defined by $J$ might be different from $E_1 (I), \ldots, E_k (I)$
but lie in the same homology classes which are independent of the complex structure.
We will denote the cohomology classes that are Poincar\'e-dual to these homology classes
by $[E_1],\ldots, [E_k]\in H^2 (\tM,\Z)$.

Now let us briefly recall the notion of a (simultaneous) symplectic blow-up at $k$ points --
for details see \cite{_McD-Polt_}, \cite{_McD-Sal-intro_}.
Denote by $B^{2n} (r)$ a round closed ball of
radius $r$ in the standard symplectic $\R^{2n}$. Given a symplectic
manifold $(M^{2n},\omega)$ and a symplectic embedding
$\iota: \bigsqcup\limits_{i=1}^k B^{2n} (r_i) \to (M,\omega)$,
one can construct a new manifold, diffeomorphic to $\tM$, by
removing $\iota (\bigsqcup\limits_{i=1}^k B^{2n} (r_i))$ from $M$ and contracting the boundary
of the resulting manifold along the fibers of the fibration induced
by $\iota$ from the Hopf fibration on the boundary of $B^{2n} (r)$.
The form $\omega$ is then extended in a certain way to a symplectic
form $\tomega$ on $\tM$ and the resulting symplectic manifold
$(\tM,\tomega)$ is independent of the extension choices, up to a
symplectic isotopy (though it does depend on $r_1,\ldots, r_k$ and $\iota$), and
is called {\bf the symplectic blow-up of $(M,\omega)$ along
$\iota$}. Alternatively, one can describe the construction of
$(\tM,\tomega)$ as follows: extend $\iota$ to a symplectic
embedding of the union of a slightly larger closed balls, remove the interior of the image of the union of the larger balls from $M$
and glue in symplectically the disjoint union of appropriate neighborhoods of $\C P^{n-1}$ in $\C P^n$ with the
standard symplectic form on them (normalized so that
its integral over the projective line is equal to $\pi$).

The cohomology class of $\tomega$ is given by
\[
[\tomega] = \Pi^* [\omega] - \pi\sum_{i=1}^k r_i^2 [E_i].
\]

McDuff and Polterovich discovered in \cite{_McD-Polt_} that, in fact, the existence of a symplectic form on
$\tM$ in the cohomology class $\Pi^* [\omega] - \pi\sum_{i=1}^k r_i^2 [E_i]$ satisfying certain additional conditions is sufficient for
the existence of a symplectic embedding $\iota: \bigsqcup\limits_{i=1}^k B^{2n} (r_i) \to (M,\omega)$.
We will state their result in the case when the symplectic manifold $M$ is K\"ahler (since we are going to
work only with such manifolds).


\hfill


\theorem (McDuff-Polterovich, \cite[Cor. 2.1.D and Rem.
2.1.E]{_McD-Polt_}) \label{_McD_P_main_Theorem_}
\\
Let $M$, $\dim_\R M
= 2n$, be a closed connected complex manifold equipped with a K\"ahler form
$\omega$. Let $k\in\N$ and let $\tM$, $\Pi^* : H^2 (M,\R)\to H^2
(\tM,\R)$ and $[E_i]\in H^2(\tM, \Z)$, $i=1,\ldots,k$, be defined as
above. Let $r_1,\ldots, r_k$ be a collection of positive numbers.
Assume there exists a complex structure $I$ of K\"ahler type on $M$
tamed by $\omega$ and a symplectic form $\tomega$ on $\tM$ taming
$\tI$ so that $[\tomega] = \Pi^* [\omega] - \pi \sum_{i=1}^k r_i^2
[E_i]$. Let $\Gamma\subset M$ be a (possibly empty) closed complex
submanifold of $(M,I)$ that does not intersect the image of
$\bigcup_{i=1}^k  E_i (I)$ under $\Pi_I$. Then $(M\setminus
\Gamma,\omega)$ admits a symplectic embedding of
$\bigsqcup\limits_{i=1}^k B^{2n} (r_i)$.
\endproof

\hfill


\remark \label{_some forms are Kahler_Remark_}
\\
The proof of
\ref{_McD_P_main_Theorem_} in \cite{_McD-Polt_} actually shows that
for any sufficiently small $c_1,\ldots, c_k >0$ the cohomology class
$\Pi^*  [\omega] - \sum_{i=1}^k \pi c_i [E_i]$ is K\"ahler (with
respect to $\tI$). Thus $\Kah (\tM, \tI)$ is non-empty for any
complex structure $I$ of K\"ahler type on $M$ or, in other words, if
$I$ is of K\"ahler type on $M$, then $\tI$ is of K\"ahler type on
$\tM$.


\hfill


\theorem\label{_Kahler_McD_Polt_Theorem_}\\
Let $(M, I, \omega)$, $\dim_\R M = 2n$, be a closed connected K\"ahler manifold. Let $k\in\N$ and let $\tM$,
$\Pi^* : H^2 (M,\R)\to H^2 (\tM,\R)$,  $[E_1], \ldots, [E_k]\in H^2(\tM, \Z)$,
$r_1,\ldots, r_k >0$ be as above.
Let $\Gamma\subset M$ be a (possibly empty)
closed complex submanifold of $(M,I)$ that does not intersect the image of $\cup_{i=1}^k  E_i (I)$ under $\Pi_I$.
Assume that there exists a complex structure $J$ of K\"ahler type on $M$ which is tamed by $\omega$ so that
$\Pi^*[\omega]_J^{1,1} - \pi \sum_{i=1}^k r_i^2 [E_i]\in \Kah (\tM, \tJ)$.

Then $(M\setminus \Gamma, \omega)$ admits a symplectic embedding of
$\bigsqcup\limits_{i=1}^k B^{2n} (r_i)$.


\hfill


\noindent {\bf Proof:}

First, let us remark that any complex
structure on $M$ defines a complex structure of K\"ahler type on
$\tM$ (see \ref{_some forms are Kahler_Remark_}) that defines a
Hodge decomposition on $H^2 (\tM,\R)$. Under the identification $H^2
(\tM,\R) = H^2 (M,\R) \oplus \Span_\R \{ [E_1],\ldots,[E_k]\}$ the
homomorphism $\Pi^*$ (which is independent of the complex structure
on $M$) acts as an identification of $H^2 (M,\R)$ with the first
summand. The classes $[E_1],\ldots,[E_k]$ are all of type $(1,1)$.
This identification preserves the Hodge types (with respect
to the complex structures on $M$ and $\tM$).

Since, by our assumption, the cohomology class $\Pi^*[\omega]_J^{1,1} - \pi \sum_{i=1}^k r_i^2 [E_i]$ lies in $\Kah (\tM, \tJ)$,
it can be represented by a K\"ahler form $\talpha$ on $(\tM, \tJ)$.

Note that $\Pi^* [\omega]_J^{1,1}\in H^2 (\tM,\R)$ is of type $(1,1)$ with respect to $\tJ$.
Hence, the class $\Pi^* [\omega] - \Pi^* [\omega]_J^{1,1}\in H^2 (\tM,\R)$ is of type $(2,0)+(0,2)$ with respect to $\tJ$
and can be represented as $\Pi^* b$ for a $(2,0)+(0,2)$-class $b\in H^2 (M,\R)$ with respect to $J$. Represent $b$ by a
closed real-valued form $\beta$ on $M$ of type $(2,0)+(0,2)$ with respect to $J$. Then
the class $\Pi^* [\omega] - \Pi^* [\omega]_J^{1,1}$ is represented by a closed real-valued form
$\Pi_J^* \beta$ on $\tM$ of type $(2,0)+(0,2)$ with respect to $\tJ$.

Set $\tomega: = \talpha + \Pi_J^* \beta$. By \ref{_tamed+2,0_Proposition_}, parts (A) and (C), the form $\tomega$ is symplectic and tames $\tJ$.
The cohomology class of $\tomega$ can be written as
\[
[\tomega] = [\talpha] + [\Pi_J^* \beta] = \Pi^* [\omega]_J^{1,1} - \pi \sum_{i=1}^k r_i^2 [E_i] + \Pi^* [\omega] - \Pi^* [\omega]_J^{1,1} =
\]
\[
=  \Pi^* [\omega] - \pi \sum_{i=1}^k r_i^2 [E_i].
\]
Now we can apply \ref{_McD_P_main_Theorem_} {\it with $J$ instead of $I$}, which yields the needed claim.
\endproof


\hfill


A necessary condition for $\Pi^*[\omega] - \pi \sum_{i=1}^k r_i^2 [E_i]$
to be K\"ahler is $$\langle (\Pi^*[\omega] - \pi \sum_{i=1}^k r_i^2 [E_i])^n,
[\tM]\rangle >0.$$ The following proposition shows that in terms of
symplectic packings the latter inequality means the following simple fact: if a finite
disjoint union of closed balls is symplectically embedded in
$(M,\omega)$, then its total volume is less than the volume of
$(M,\omega)$.


\hfill

\proposition\label{_highest power_Proposition_}\\
With the notation as in \ref{_McD_P_main_Theorem_},
\[
\langle (\Pi^* [\omega] - \pi \sum_{i=1}^k r_i^2 [E_i])^n, [\tM]\rangle = \int_M \omega^n -  \pi^n \sum_{i=1}^k r_i^{2n} =
\]
\[
= \int_M \omega^n - \Vol (\bigsqcup\limits_{i=1}^k B^{2n} (r_i)).
\]


\hfill


\noindent {\bf Proof:}

Note that for all $i=1,\ldots,k$ we have $\Pi^*
[\omega] \cup [E_i] = 0$, as well as $[E_i]^n = -1$, if $n$ is even, and
$[E_i]^n = 1$, if $n$ is odd. Also note that $[E_i]\cup [E_j] = 0$
for all $i\neq j$. Finally, recall that the symplectic volume of
$B^{2n} (r)$ equals $\pi^n r^{2n}$. The claim follows directly from
these observations.
\endproof


\hfill

\subsection{Demailly-Paun theorem and the K\"ahler cone}
\label{_Demailly_Paun_Section_}

Our results depend crucially on the following deep result by Demailly and Paun.


\hfill

\theorem (Demailly-Paun, \cite{_Demailly_Paun_})
\label{_Dem_Paun_cone_Theorem_}\\
Let $N$ be a closed connected K\"ahler manifold. Let $\hat K(N)\subset H^{1,1}(N,\R)$ be
a subset consisting of all (1,1)-classes $\eta$ which satisfy
$\langle \eta^m, [Z]\rangle >0$ for any homology class $[Z]$ realized by a complex subvariety
$Z\subset N$ of complex dimension $m$. Then the K\"ahler cone of $N$ is one of the
connected components of $\hat K(N)$. \endproof


\hfill


For Campana simple manifolds this theorem can be used to study
the K\"ahler cone of a blow-up.


\hfill

\theorem \label{_Kah_cone_blow_up_simple_Theorem_}\\
Let $(M,I)$, $\dim_\C M = n$, be a Campana simple  closed connected K\"ahler manifold.
Consider a complex blow-up $\tM$ of $(M,I)$ at $k$ Campana-generic points $x_1,\ldots, x_k$.
Define $\Pi_I: \tM\to M$, $E_i := \Pi_I^{-1} (x_i)$ and $[E_i]\in H^2 (\tM,\Z)$, $i=1,\ldots,k$,
as above.

Assume that $\eta$ is a K\"ahler class in $H^2 (M,\R)$.
Then, given $c_1,\ldots,c_k\in \R$,
the following claims are equivalent:
\begin{description}
\item[(A)] The cohomology class $\teta:=\Pi^* \eta - \sum_{i=1}^k c_i [E_i] \in H^2(\tM,\R)$ is K\"ahler.
\hfill
\item[(B)] The conditions (B1) and (B2) below are satisfied:\\
(B1) All $c_i$ are positive. \\
(B2) $\langle \teta^n, [\tM]\rangle >0$.
\end{description}


\hfill


\noindent
{\bf  Proof of (A) $\Rightarrow$ (B).}\\
The implication (A) $\Rightarrow$ (B2) is obvious. To prove (A) $\Rightarrow$ (B1) note that, since $\teta$ is K\"ahler,
for each $i=1,\ldots,k$ we have
$$0<\int_{E_i (I)} \teta^{n-1} = \int_{E_i (I)} (-c_i [E_i])^{n-1},$$
and since the restriction of $-[E_i]$ to $E_i (I)$ is a positive multiple of the Fubini-Study form on the exceptional divisor
(see the discussion on the symplectic blow-up in Section~\ref{_blow_ups_McDuff_Polterovich_Section_}) and the integral
of the exterior power of the latter form over $E_i (I)$ is positive, we readily get
 that $c_i >0$.
\endproof


\hfill


\noindent {\bf Proof of (B) $\Rightarrow$ (A).}

Assume (B1) and (B2) are satisfied.

Since $x_i$ are Campana-generic,
any connected
proper complex subvariety of $(\tM, \tI)$ is either contained
in an exceptional divisor $E_i$, or does not intersect any exceptional divisor.

Since $\eta\in H^2 (M,\R)$ is a K\"ahler cohomology class with respect to $I$, we have
$\langle \teta^m, Z\rangle = \langle \eta^m, \Pi_I (Z)\rangle >0$ for any complex subvariety $Z\subset (\tM, \tI)$, $\dim_\C Z = m$, that does not intersect the exceptional divisors.

On the other hand, note that for each $i=1,\ldots,k$ the restriction of the cohomology class $[E_i]$ to the submanifold $E_i (I)$ is a positive multiple of $- [\omega_{E_i (I)}]$,
where $\omega_{E_i (I)}$ is the Fubini-Study form, and therefore
(B1) yields that $\langle \teta^m, [Z]\rangle >0$ for any complex variety $Z\subset (\tM,\tI)$,
$\dim_\C Z = m$, lying in $E_i (I)$.

Thus $\langle \teta^m, Z\rangle >0$ for any complex subvariety $Z\subsetneq (\tM,\tI)$, $\dim_\C Z = m$.
This shows that for any positive $c_1,\ldots, c_k$ the class $\teta=\eta - \sum_{i=1}^k c_i [E_i]$ lies in
$\hat K(\tM)$ as long as it satisfies (B2).

By \ref{_Dem_Paun_cone_Theorem_}, in order to show that $\teta$ is K\"ahler, it remains to check that there exists a K\"ahler form
in the connected component of $\hat K(\tM)$ containing $\teta$.

Indeed, similarly to \ref{_highest power_Proposition_}, one gets that (B2) is equivalent to the condition
\[
\sum_{i=1}^k c_i^n < \langle \teta^n, [\tM]\rangle.
\]
If this condition holds for $c_1,\ldots, c_k >0$, it also holds for
$\epsilon c_1,\ldots, \epsilon c_k$ for any $\epsilon\in (0,1]$. The
numbers $\epsilon c_1,\ldots, \epsilon c_k$ are still positive and
therefore, by the argument above, for any $\epsilon\in (0,1]$ the
class $\teta_\epsilon :=\Pi^* \eta - \epsilon \sum_{i=1}^k c_i
[E_i]$ also lies in $\hat K(\tM)$. But, as it is explained in
\ref{_some forms are Kahler_Remark_}, for any sufficiently small
positive $\epsilon$ the class $\teta_\epsilon$ is K\"ahler. Thus,
$\teta$ lies in the same connected component of $\hat K(\tM)$ as a
K\"ahler class $\teta_\epsilon$. Therefore, by
\ref{_Dem_Paun_cone_Theorem_}, $\teta$ is K\"ahler.
\endproof


\hfill


\subsection{Unobstructed symplectic packings by balls for Campana simple complex manifolds
and their limits}
\label{_full_packing_Campana_simple_Section_}


The following result is a slightly stronger version of
\ref{_Campana_full_packing_Theorem_}.


\hfill


\theorem\label{_Campana_simple_FSP_precise_version_Theorem_}\\
Let $(M,I,\omega)$ be a closed connected K\"ahler manifold. Assume
$I$ admits a smooth local deformation $\{ I_t\}$, $t\in U$,
$I=I_{t_0}$, and there is a sequence $\{ t_l\}_{l\in\N} \to t_0$ in
$U$ such that $I_{t_l}$ is Campana simple for all $l\in \N$.
Let $\Gamma\subset (M,I)$ be a (possibly empty)
closed complex submanifold.

Then the symplectic packings of $(M\setminus\Gamma, \omega)$ by balls are unobstructed.


\hfill


\noindent {\bf Proof or \ref{_Campana_simple_FSP_precise_version_Theorem_}:}

Consider a disjoint union $\bigsqcup\limits_{i=1}^k B^{2n} (r_i)$ whose total symplectic volume is less than the symplectic volume of $M$, that is,
\begin{equation}
\label{eqn-total-vol-balls-less-than-vol-M}
\pi^n \sum_{i=1}^k r_i^{2n} < \Vol (M) = \langle [\omega]^n, [M]\rangle.
\end{equation}
We need to show that it admits a symplectic embedding into $(M\setminus\Gamma,\omega)$.

It follows from \ref{_Kodaira_stabi_Theorem_} and the assumptions of the theorem that for
a sufficiently large $l\in \N$
and the corresponding Campana simple complex structure
$J:=I_{t_l}$ the cohomology class $[\omega]_J^{1,1}$ is K\"ahler (with respect to $J$).
Note that $J$ can be chosen arbitrarily $C^\infty$-close to $I$.
In particular, we can assume that
$J$ is
tamed by $\omega$ (since $I$ is tamed by $\omega$)
and that
\begin{equation}
\label{eqn-total-vol-balls-less-than-vol-M-J}
\pi^n \sum_{i=1}^k r_i^{2n} < \langle ([\omega]_J^{1,1})^n, [M]\rangle,
\end{equation}
because of \eqref{eqn-total-vol-balls-less-than-vol-M} and \ref{_Kodaira_stabi_Theorem_}, combined with the fact that $[\omega]=[\omega]_I^{1,1}$.

Choose $k$ Campana-generic points $x_1,\ldots,x_k \in (M,J)$ lying
in $M\setminus \Gamma$ and consider the complex blow-up $(\tM, \tJ)$
of $(M,J)$ at those points. By
\ref{_Kah_cone_blow_up_simple_Theorem_} applied to the K\"ahler
class $[\omega]_J^{1,1}$, the cohomology class $\Pi^*
[\omega]_J^{1,1} - \pi \sum_{i=1}^k r_i^2 [E_i]$ is K\"ahler with respect
to $\tJ$ (note that, by \ref{_highest power_Proposition_}, the
condition (B2) in \ref{_Kah_cone_blow_up_simple_Theorem_} is
equivalent to \eqref{eqn-total-vol-balls-less-than-vol-M-J}).
Therefore, by \ref{_Kahler_McD_Polt_Theorem_}, $(M\setminus \Gamma,
\omega)$ admits a symplectic embedding of $\bigsqcup\limits_{i=1}^k
B^{2n} (r_i)$.
\endproof


\section{Symplectic packing by arbitrary shapes}
\label{_Packing_by_arb_shapes_main_Section_}


Let $M$ be either an oriented torus $T^{2n}$, $n\geq 2$, or, respectively, a closed connected oriented manifold admitting IHS
hyperk\"ahler structures compatible with the orientation. Without loss of
generality we are going to prove the results for symplectic forms on
$M$ of total volume 1.

\subsection{Semicontinuity of symplectic packing constants in families of sym\-plec\-tic forms}
\label{_semicont_of_sympl_packing_section_}

 Let $\cF$ denote the space of K\"ahler forms on $T^{2n}$ (respectively, hyperk\"ahler
forms on $M$ in the hyperk\"ahler case) of total volume $1$. Equip
$\cF$ with the $C^\infty$-topology. The group $\Diff^+$ of
orientation-preserving diffeomorphisms of $M$ acts on $\cF$ and the
function $\omega\mapsto \nu (M,\omega, V)$, defined in Subsection
\ref{_arbi_shapes_Subsection_}, is clearly invariant under the
action.


\hfill


\proposition
\label{_semicontinuity_Proposition_}\\
The function $\omega\mapsto \nu (M,\omega, V)$ on $\cF$ is lower
semicontinuous.


\hfill

\theorem
\label{_ergodicity_Theorem_}\\
Let $\omega\in \cF$ be a symplectic form such that the cohomology
class $[\omega]$ is not proportional to a rational one. Then the orbit of $\omega$ under
action of $\Diff^+$ is dense in $\cF$ in the toric case and in the
connected component $\cF^0$ of $\cF$ containing $\omega$ in the
hyperk\"ahler case.


\hfill


We will prove \ref{_ergodicity_Theorem_} in
Subsection~\ref{_Ergodic_Subsection_}.


\hfill


\noindent
{\bf Proof of \ref{_Packing_by_arb_shapes_main_Theorem_}.}

Since the orbit
of $\omega_1$ is dense in $\cF$ and $\nu (M,\cdot, V)$ is constant
on it, we get, by the lower semicontinuity of $\nu (M,\cdot, V)$,
that
$$\nu (M,\omega_2, V) \leq \nu (M,\omega_1, V).$$
Switching $\omega_1$ and $\omega_2$ and applying the same argument
we get
$$\nu (M,\omega_1, V) \geq \nu (M,\omega_2, V).$$
Thus
$$\nu (M,\omega_1, V) = \nu (M,\omega_2, V).$$
\endproof


\hfill
\eject


\noindent
{\bf Proof of \ref{_semicontinuity_Proposition_}.}

Let us start with a number of preparations.

Fix Riemannian metrics on $U$ and on $M$.
These Riemannian metrics induce $C^0$-norms on the spaces of vector
fields and differential forms defined on an open neighborhood of $V$
in $U$ and on $M$. Abusing the notation, we will denote all these
norms by the same symbol $||\cdot||$:
$$||v||:=\max_x |v(x)|$$
for a vector field $v$ and
$$||\Omega|| := \max_{ ||v_1||,\ldots,||v_l||\leq 1}  |\Omega (v_1,\ldots,v_l)|$$
for a differential $l$-form $\Omega$ (for vector fields and forms on
a neighborhood of $V$ we take the maximums only over $V$ -- recall
that $V$ is compact).


\hfill


\lemma
\label{_small_primitive_for_a_small_exact_form_Lemma_}\\
Let $U$, $\dim_\R U = 2n$, be an open,
possibly disconnected
symplectic manifold, and
let $V\subset U$, $\dim_\R V = 2n$, be a compact
submanifold with a
piecewise smooth boundary.
Given an exact 2-form $\Omega$ on a neighborhood of $V$, one can
choose a 1-form $\sigma$ on the same neighborhood so that
$d\sigma=\Omega$ and $\| \sigma\| \leq C_1 \|\Omega\|$
for some constant $C_1>0$ depending only on $V$.

\hfill

\noindent
{\bf Proof:}

We present only an outline of the proof leaving the technical details to the reader.

We triangulate $V$ and proceed by induction on the number of
$2n$-dimensional simplices. In the case of one simplex the claim follows from an
explicit formula for the primitive of an exact form on a star-shaped
domain in a Euclidean space appearing in the proof of the classical
Poincare lemma (see e.g. \cite{_Spivak_}).

Assume now that the
result holds for any manifold whose triangulation consists of $k$
simplices and consider a manifold which is the union of $k+1$
simplices. Apply the induction assumption to the union $A$ of the
first $k$ simplices and, separately, to the $(k+1)$-st simplex $B$.

If $A\cap B=\emptyset$ the claim is obvious. Therefore we may assume that $A\cap B\neq \emptyset$.
We get two small 1-forms $\sigma_1$ and $\sigma_2$ defined on open
neighborhoods $Z$ and $W$ of $A$ and $B$ so that $d\sigma_1=\Omega$
and $d\sigma_2=\Omega$. Thus on $Z\cap W$ the 1-form
$\sigma_1-\sigma_2$ is exact and small. Without loss of generality,
we may assume that $W$ is a ball and $Z\cap W$ is, topologically,
either a ball or a spherical shell. In either case it is not hard to
see that $\sigma_1-\sigma_2$ can be written on $Z\cap W$ as
$\sigma_1-\sigma_2=dh$ for a $C^1$-small function $h$. Extend the
function $h$ from $Z\cap W$ to $W$ keeping it $C^1$-small (this is
not hard to do, since $Z\cap W$ is a ball or a spherical shell
inside the ball $W$). This allows to extend $\sigma_1-\sigma_2$ to a
small exact 1-form on $W$. Thus $\sigma_1$ (which is equal to
$\sigma_1= \sigma_2 - (\sigma_2-\sigma_1)$ on $Z\cap W$) can be
extended to a small $1$-form on the open set $Z\cup W$ (which is a
neighborhood of our original manifold $V$) so that $d\sigma_1
=\Omega$ everywhere.
The loss of ``smallness" of the differential forms at each
step above is by a factor that depends only on the
geometry of $V$ and not on the differential forms. Setting
$\sigma:=\sigma_1$ finishes the proof.
\endproof


\hfill


Having fixed $V$ and $\eta$, denote for brevity
$\nu(\omega):= \nu (M,\omega, V)$. Consider an arbitrary
$\omega_0\in \cF$. Let us prove that $\nu$ is lower semicontinuous
at $\omega_0$. If $\nu (\omega_0) = -\infty$, the claim is obvious,
so we can assume without loss of generality that $\nu (\omega_0)
\neq -\infty$, meaning that there exist symplectic embeddings $(V,
\alpha\eta)\to (M,\omega_0)$ for some $\alpha >0$. To prove the
lower semicontinuity of $\nu$ at $\omega_0$ it suffices to prove the
following claim:

For any symplectic embedding $f: (\cU, \alpha\eta)\to (M,\omega_0)$,
where $\cU\subset U$ is an open neighborhood of $V$ in $U$, and any
sufficiently small $\varepsilon >0$ there exists $\delta = \delta
(f, \omega_0,\varepsilon) >0$ such that for any $\omega\in \cF$, $|| \omega -
\omega_0||\leq \delta$, the following two conditions are satisfied:

\medskip
\noindent (A) $\displaystyle \frac{\Vol (V,\alpha \eta)}{\Vol
(M,\omega)} > \frac{\Vol (V,\alpha\eta)}{\Vol (M,\omega_0)} -
\varepsilon$.

\medskip
\noindent (B) There exists a symplectic embedding $g: (\cU',
\alpha\eta)\to (M,\omega)$, where $\cU'\subset \cU$ is a possibly
smaller neighborhood of $V$.
\medskip

In order to prove the claim let us fix $\varepsilon >0$ and choose
$\delta>0$ (depending on $\omega_0$) so that (A) is satisfied for all  $\omega\in \cF$, $|| \omega -
\omega_0||\leq \delta$ -- this is, of course, not a
problem. Consider now a form $\omega\in
\cF$, $|| \omega - \omega_0||\leq \delta$, and let us show that (B) holds for $\omega$: namely, we will construct $g$ by the
classical Moser method \cite{_Moser_} as follows. If $\delta$ is
sufficiently small (depending on $\omega_0$ and $f$), then, since
$f^* \omega_0 = \alpha \eta$ is symplectic, the straight path
$\theta_t := f^* \omega_0 + t f^* (\omega - \omega_0)$, $0\leq t\leq
1$, connecting $f^* \omega_0$ and $f^* \omega$ is formed by
symplectic forms on $V$. Since $H^2 (V,\R)=0$, the form $f^* (\omega
- \omega_0)$ is an exact. By
\ref{_small_primitive_for_a_small_exact_form_Lemma_}, one can choose
a 1-form $\sigma$ on $\cU$ so that $f^* (\omega - \omega_0) =
d\sigma$ and $\| \sigma\| \leq C_2 \|\omega-\omega_0\|$ for some
constant $C_2>0$ depending only on $V$, $\omega_0$ and $f$. Let
$v_t$ be the vector field on $\cU$ defined by $\theta_t (v_t, \cdot)
= \sigma (\cdot)$. Then $\max_{0\leq t\leq 1} ||v_t || \leq C_3
||\omega-\omega_0||\leq C_3\delta$ for some constant $C_3>0$
depending only on $V$, $\omega_0$ and $f$. Therefore if $\delta$ is
sufficiently small, the time-$[0,1]$ flow of $v_t$ yields a
well-defined map $\psi: \cU'\to \cU$ for some smaller neighborhood
$\cU'\subset \cU$ of $V$. Moser's argument \cite{_Moser_} shows that
$\psi^* (f^* \omega) = f^* \omega_0$. This implies that $g:=f\circ
\psi : (\cU', \alpha \eta)\to (M,\omega)$ is a symplectic embedding:
indeed,
\[ g^* \omega = (f\circ \psi)^* \omega = \psi^* (f^* \omega) = f^* \omega_0 = \alpha \eta. \]
\endproof


\hfill


\noindent
{\bf Proof of \ref{_full_packing_of_tori_by_polydisks_Corollary_}.}

Let $(V,\eta)$ be the disjoint union of $k$ copies of $(B^{2n_1}
(R_1)\times \ldots B^{2n_l} (R_l),dp\wedge dq)$, where $dp\wedge dq$
is the standard symplectic form on $\R^{2n}$. For brevity set $\nu
(\omega):= \nu (T^{2n},\omega, V)$. We need to show that $\nu
(\omega) = 1$.

For each $i=1,\ldots,l$ set $v_i:= \Vol (B^{2n_i} (R_i),
\Omega_{n_i})$, where $\Omega_{n_i}$ is the standard symplectic form
on $\R^{2n_i}$. Note that
$$v_i = \Vol (B^{2n_i} (R_i), \Omega_{n_i})
= \pi^{n_i} n_i! R_i^{2n_i},$$
$$dp\wedge dq =
\Omega_{n_1}\oplus\ldots\oplus \Omega_{n_l}$$
and
$$\Vol (V,\eta) = k\Vol (B^{2n_1} (R_1)\times \ldots B^{2n_l} (R_l),dp\wedge dq)
= kN\prod_{i=1}^l v_i,$$ where
$$N:=\frac{n!}{n_1!\cdot\ldots\cdot
n_l!}.$$

Assume without loss of generality that
$$\Vol (T^{2n},\omega) = k\Vol (B^{2n_1} (R_1)\times \ldots B^{2n_l} (R_l),dp\wedge dq) = kN v_1\cdot \ldots \cdot v_l
=1.$$

Given $m\in\N$ and $w_1,\ldots, w_m >0$, set $\bw:= (w_1,\ldots,
w_m) \in \R^m$, $v_\bw := m! w_1\cdot\ldots \cdot w_m$, and denote
by $\omega_\bw$ the symplectic form $\omega_\bw = \sum_{i=1}^m w_i
dp_i\wedge dq_i$ on the torus $T^{2m} = \R^{2m}/\Z^{2m}$. Note that
the form $\omega_\bw$ is K\"ahler (since it is linear) and $\Vol
(T^{2m}, \omega_\bw) = v_\bw$.

Set $k_1:=k$, $k_2=\ldots=k_l:=1$, so that
$$k = k_1\cdot\ldots\cdot k_l.$$

Given $0< \alpha\leq 1$, one can choose $\bw_i\in \R^{2n_i}$,
$i=1,\ldots,l$, depending on $\alpha$ (for brevity we suppress this
dependence in the notation below) so that the following conditions
hold:

\smallskip
\noindent (A) $\prod_{i=1}^l v_{{\bw_i}} = \prod_{i=1}^l k_i v_i =
1/N$.

\smallskip
\noindent (B) $v_{{\bw_i}} > \alpha k_i v_i$ for all $i=1,\ldots,l$.

\smallskip
\noindent (C) The vector $\bbw := (\bw_1,\ldots,\bw_k) \in \R^{2n}$ is
not proportional to a vector with rational coordinates.

\smallskip
Condition (C) can be achieved since the set of vectors that are not
proportional to a vector with rational coordinates is dense in the
set
$$\{ (w_1,\ldots,w_{2n})\in\R^{2n}\ |\ w_1,\ldots,w_{2n}>0,
w_1\cdot\ldots \cdot w_{2n} = C\}$$
for any $C>0$.

Consider the symplectic form $\omega_{\bbw}$ on $T^{2n}$ -- it is K\"ahler (since it is a linear symplectic form)
and its cohomology class is not proportional to a rational one. Note that, by condition (A),
\[
\int_{T^{2n}} \omega_{\bbw}^n = N \prod_{i=1}^l v_{{\bw_i}}  = \int_{T^{2n}} \omega^n = 1.
\]

Thus to prove that $\nu (\omega) = 1$ it is enough to show that $\nu
(\omega_{\bbw}) \to 1$ as $\alpha\to 1$. Indeed, by
\ref{_Packing_by_arb_shapes_main_Theorem_} and in view of condition
(C), $\nu (\omega_{\bbw})$ is constant as a function of $\alpha$ for all $0< \alpha\leq 1$ (recall
that $\bbw$ depends on $\alpha$) and equal to $\nu (\omega)$.

Now note that, by condition (B), for any $0< \alpha\leq 1$ for all $i$
$$\Vol (T^{2n_i},\omega_{{\bw_i}}) = v_{{\bw_i}}
> \alpha k_i v_i  = \Vol \bigg( \bigsqcup_{k_i} (B^{2n_i} (R_i), \alpha \Omega_{n_i}) \bigg),$$
where $\bigsqcup_{k_i}$ denotes the disjoint union of $k_i$ copies
of the ball. Therefore, by \ref{_FSP_main_Theorem_}, there exists a
symplectic embedding $$f_i: \bigsqcup_{k_i} \big(B^{2n_i} (R_i),
\alpha \Omega_{n_i}\big) \to (T^{2n_i}, \omega_{\bw_i}).$$
Accordingly, the direct product of all such embeddings $f_i$,
$i=1,\ldots,l$, is a symplectic embedding
$$f: (V,\alpha\eta) \to
(T^{2n},\omega_{\bbw})$$ of $k_1\cdot\ldots\cdot k_l=k$ disjoint
equal copies of the polydisk $(B^{2n_1} (R_1)\times \ldots \times
B^{2n_l} (R_l),dp\wedge dq)$ into $(T^{2n},\omega_{\bbw})$. The
fraction of the volume of $(T^{2n},\omega_{\bbw})$ filled by the
image of $f$ tends to $1$ as $\alpha\to 1$. In other words, $\nu
(\omega_{\bbw})$ converges to $1$ as $\alpha\to 1$ which yields the
needed result. \endproof

\subsection{Ergodic action on the symplectic Teichm\"uller space -- the proof of \ref{_ergodicity_Theorem_}}
\label{_Ergodic_Subsection_}

The proof of \ref{_ergodicity_Theorem_} follows the
same lines as the proof of the ergodicity theorem in
\cite{_Verbitsky:ergodic_} and \cite{_Verbitsky:ICM_}.

First, let us make a few preparations.

Assume $G$ is a Lie group and $\Gamma\subset G$ is a discrete subgroup. The restriction of the right Haar measure on $G$ to a fundamental domain of the right action of $\Gamma$ on $G$ induces a measure on the space $G/\Gamma$. The subgroup $\Gamma$ is called a {\bf lattice} if the total measure of $G/\Gamma$ is finite.

Assume now that $G\subset SL (m,\R)$ -- that is, $G$ is a linear Lie group.

A {\bf $\Q$-character} on $G$ is a homomorphism $G\to
 \R_{>0}$ which is defined by algebraic equations with rational
 coefficients on the entries of the real $m\times
 m$-matrices in $G$.

The group $G$ is called an {\bf algebraic group}
(respectively, an {\bf algebraic $\Q$-group})
if it is defined by algebraic equations with real
(respectively, rational)
coefficients on the entries of the
matrices.

A lattice
$\Gamma\subset G$ is called {\bf arithmetic} if $\Gamma \cap
SL (m,\Z)$ has a finite index in $\Gamma$ and in $SL (m,\Z)$.

The following theorem belongs to A.Borel and Harish-Chandra.


\hfill


\claim (Borel-Harish-Chandra theorem, \cite[Thm. 9.4]{_Borel-HC1962_})
\label{_Borel-Harish-Chandra_Claim_}

Let $G\subset SL (m,\R)$ be an algebraic $\Q$-group. Then
$G\cap SL (m,\Z)$ is a Lie lattice in $G$ (or, equivalently, $G$ admits an arithmetic lattice)
if and only if $G$ admits no non-trivial $\Q$-characters.
\endproof


\hfill


\proposition
\label{_lattice-identity-component_Proposition_}

Let $G'\subset SL (m,\R)$ be an algebraic $\Q$-group and $G$ its identity component (in the Lie group topology). Assume
$G$ does not admit non-trivial $\Q$-characters.

Then $G\cap SL (m,\Z)$ is a Lie lattice in $G$.


\hfill


\noindent
{\bf Proof of \ref{_lattice-identity-component_Proposition_}:}

Since $G$ does not admit non-trivial $\Q$-characters, $G'$ does not admit them either. Therefore, by Borel and Harish-Chandra theorem
(\ref{_Borel-Harish-Chandra_Claim_}), $G'_\Lambda$ is a Lie lattice in $G'$.

The identity component $G$ is a normal subgroup of $G'$ and its index is finite, since, by Whitney's theorem \cite{_Whitney_}, the real algebraic
subvariety of a Euclidean space has finitely many connected components. This implies that $G\cap SL (m,\Z)$ has a finite index in $G'\cap SL (m,\Z)$.
Since the indices $|G':G|$ and $|(G'\cap SL (m,\Z)):(G\cap SL (m,\Z))|$ are finite and $G'\cap SL (m,\Z)$ is a Lie lattice in $G'$, one easily gets that $G\cap SL (m,\Z)$ is a Lie lattice in $G$.
\endproof


\hfill

\definition
\\
Let $G$ be a real Lie group. We say that $g\in G$ is {\bf
unipotent}, if $g=e^h$ for a nilpotent element $h$ the Lie algebra
of $G$. A group $G$ is said to be {\bf generated by unipotents}, if
it is multiplicatively generated by unipotent elements.

A {\bf unipotent one-parametric subgroup of $G$} is a subgroup of the form $\{ e^{th}\}_{t\in\R}$ for an $ad$-nilpotent element $h$ of the Lie algebra of $G$.


\hfill


\proposition
\label{_SL-n-R-SO-p-q-generated-by-unipotents_Proposition_}

The groups $SL (m,\R)$, $n\geq 2$, and $SO^+ (p,q)$
(the identity component of $SO(p,q)$), $p,q\in\Z_{>0}$, $p+q>2$,  are generated by their algebraic one-parametric unipotent subgroups.


\hfill


\noindent
{\bf Proof of \ref{_SL-n-R-SO-p-q-generated-by-unipotents_Proposition_}:}

The group $SL(2,\R)$, hence $SL(m,\R)$ for any $m\geq 2$,
contains a non-trivial algebraic one-parametric uni\-potent subgroup
\[
U:=\left\{
\begin{pmatrix}
1 & t \\ 0 & 1
\end{pmatrix},
t\in\R
\right\}.
\]
There is an isomorphism between $SL(2,\R)/\{\pm 1\}$ and $SO^+(1,2)$ defined by the following homomorphism $\mu: SL(2,\R)\to SO^+(1,2)$
with the kernel $\{\pm 1\}$: identify $\R^3$ with the space $V$ of real trace-free $2\times 2$-matrices equipped with the bilinear symmetric form $\langle A,B\rangle\to tr (AB)$, $A,B\in V$, of signature $(1,2)$; for each $C\in SL(2,\R)$ define $\mu(C)$ as the automorphism of $V$ given by $A\mapsto CAC^{-1}$, $A\in V$. The homomorphism $\mu$ is a polynomial map from $SL (2,\R)$ to $SO (1,2)$. Hence, $\mu (U)$ is a non-trivial algebraic one-parametric unipotent subgroup of $SO^+(1,2)$. This implies that for all $p,q\in\Z_{>0}$ and $p+q>2$ the group $SO^+ (p,q)$ also contains non-trivial algebraic one-parametric unipotent subgroups.

Any conjugate of an algebraic one-parametric unipotent subgroup of a real Lie group is again an algebraic one-parametric unipotent subgroup of the same Lie group. Therefore
the subgroups of $SO^+ (p,q)$ and of $SL(m,\R)$ generated by their algebraic one-parametric unipotent subgroups are normal and non-discrete.
Unless $p=q=2$, the group $SO^+ (p,q)$ is a simple Lie group and so is $SL(m,\R)$. Therefore these groups do not contain non-discrete proper normal subgroups \cite{_Ragozin_PAMS72_}. Consequently, they are generated by their algebraic one-parametric unipotent subgroups.

Let us consider the remaining case $p=q=2$. We have
$SO^+ (2,2) = \big(SL (2,\R)\times SL (2,\R)\big)/\{\pm 1\}$.
If $U$ is an arbitrary algebraic one-para\-met\-ric uni\-po\-tent subgroup of $SL (2,\R)$, then
$Id\times U$ and $U\times Id$ are algebraic one-para\-met\-ric unipotent subgroups of $SO^+ (2,2) = \big(SL (2,\R)\times SL (2,\R)\big)/\{\pm 1\}$.
Consequently, since $SL(2,\R)$ is generated by its algebraic one-parametric unipotent subgroups, so is $SO^+ (2,2)$.

This shows that the groups
$SL (m,\R)$, $m\geq 2$, and $SO^+ (p,q)$, $p,q\in\Z_{>0}$, $p+q>2$,  are generated by their algebraic one-parametric unipotent subgroups.
\endproof


\hfill


Our proof
of \ref{_ergodicity_Theorem_}
will be based on the following fundamental theorem of Ratner (see e.g. \cite{_Morris:Ratner_} for a friendly introduction to the subject).


\hfill

\theorem (Ratner's orbit closure theorem, \cite{_Ratner_Duke_1991_})
\label{_Ratner_orbit_closure_Theorem_}
\\
Let $G$ be a connected Lie group,  $H\subset G$ its subgroup
generated by unipotents and $\Gamma\subset G$ a lattice (that is, a discrete subgroup of finite covolume).
Then for any $g\in G$ one has
\[
\overline{\Gamma g H} = \Gamma g S
\]
for some closed Lie subgroup $S$, $H\subset S\subset G$.

In particular, if $H$ is a closed Lie subgroup, then the closure of
the orbit $\Gamma \cdot g H$ of $gH$ in
$G/H$ is $\Gamma (gSg^{-1}) \cdot gH$.
\endproof

\hfill

A combination of Ratner's \ref{_Ratner_orbit_closure_Theorem_} with
a result of Shah \cite[Proposition 3.2]{_Shah:uniformly_}
(cf. see \cite[Proposition 3.3.7]{_Kleinbock_etc:Handbook_})
yields a
more precise description of the group $S$ from
\ref{_Ratner_orbit_closure_Theorem_} in the case where $G$ is a
linear algebraic group and $\Gamma\subset G$ is an arithmetic
lattice.
We will state the result for $g=e$, since this is exactly what we
are going to use in our proof.


\hfill


\claim \label{_Ratner's_arithmetic_Claim_}
(\cite[Proposition 3.2]{_Shah:uniformly_},
cf. \cite[Proposition 3.3.7]{_Kleinbock_etc:Handbook_})

Let $G$ be
the identity component of
a linear algebraic $\Q$-group
$G'$, and
$\Gamma\subset G$ an arithmetic Lie lattice. Let $H\subset G$ be a closed Lie
subgroup generated by
algebraic unipotent one-parameter subgroups of $G$ contained in $H$.
Let $x:=eH\in G/H$,
where $e\in G$ is the identity of $G$.

Then the closure of the orbit $\Gamma \cdot x$ in $G/H$
is $\Gamma S\cdot x$, where $S\subset G$ is
the identity component of
the smallest
algebraic $\Q$-subgroup of
$G'$
containing $H$.
\endproof


\hfill


Let $\cF$ be the space of K\"ahler (respectively, hyperk\"ahler)
forms on $M$ defined as in
Section~\ref{_semicont_of_sympl_packing_section_} and viewed as an
infinite-dimensional Fr\'echet manifold. Let $\cF^0$ be the
connected component of $\cF$ containing the form $\omega$ as in
\ref{_ergodicity_Theorem_}. The quotient topological space
$\Teich_s:= \cF/\Diff_0$ is called {\bf the Teichm\"uller space of
symplectic structures.} Set $\Teich_s^0 := \cF^0/\Diff_0$ -- it is a
connected component of $\Teich_s$.

Let
 $\Per: \Teich_s \arrow H^2(M,\R)$ be {\bf the period map} associating
to a symplectic structure its cohomology class.
Using Moser's stability theorem for
symplectic structures, it is not hard to obtain
that $\Teich_s$ is a finite-dimensional manifold
and $\Per$
is locally a diffeomorphism \cite{_Fricke_Habermann:symp-moduli_}.

In the case when $M$ is a torus it is easy to see that $\Per$ is
a surjective local diffeomorphism\footnote{In fact, one can show that, in the case of $M=T^{2n}$, the restriction of $\Per$ to each connected component of $\Teich_s$ is a diffeomorphism. We do not know how many connected components $\Teich_s$ has.}
whose image is the set $\Theta_t \subset H^2 (T^{2n},\R)$
of cohomology classes $\eta$ such that $\int_{T^{2n}} \eta^n = 1$.
(Indeed, any such $\eta$ is the cohomology class of a linear symplectic form of total volume 1).

In the hyperk\"ahler case we will use the following theorem proved
in \cite{_Americ-Verbitsky:Teich_s_}.


\hfill

\theorem
\label{_Amerik-Verbitsky-per-map_Theorem_}\\
Let $M$, $\dim_\R M = 2n$, be a closed connected oriented manifold admitting IHS hyperk\"ahler structures (compatible with the orientation).
Then $\Per$ is an open embedding on each connected component
of $\Teich_s$ and its image is the set
$\Theta_h:=\{\ \eta \in H^2(M, \R)\ \ |\ \ q(\eta, \eta)>0,\ \int_M \eta^n = 1\ \}$,
where $q$ is the Bogomolov-Beauville-Fujiki  form defined in Section~\ref{_BBF_form_Section_}.
Moreover, $\Teich_s$ has finitely many connected components.
\endproof


\hfill


Let us return to the setup where $M$, $\dim_\R M = 2n\geq 4$, is
either an oriented torus or a closed connected oriented manifold
admitting IHS hyperk\"ahler structures (compatible with the
orientation). Let $P\subset H^2(M,\R)$
be the image of $\Per$, that is, $P=\Theta_t$ in the
torus case and $P=\Theta_h$ in the hyperk\"ahler case. The group
$\Diff^+/\Diff_0$ acts in an obvious way on $\Teich_s$ and on $H^2
(M,\R)$. The period map
 $\Per: \Teich_s \arrow H^2(M,\R)$ respects the actions.
In particular, the action of $\Diff^+/\Diff_0$ on $H^2 (M,\R)$ preserves $P$.
Let $\bGamma:= \Diff^+/\Diff_0$ in the torus case and let $\bGamma\subset \Diff^+/\Diff_0$ be
the subgroup fixing the connected component $\Teich_s^0$ (as a set) in the hyperk\"ahler case:
$\bGamma \cdot \Teich_s^0 =\Teich_s^0$.


\hfill
\bigskip

\theorem\label{_Ratner_for_H^2_Theorem_}

For any $\eta\in P$ the orbit $\bGamma\cdot \eta$ is dense in $P$ if
and only if the cohomology class $\eta$ is not proportional to a rational one.


\hfill


Before proving  \ref{_Ratner_for_H^2_Theorem_} let us see how it implies \ref{_ergodicity_Theorem_}.


\hfill


\noindent {\bf Proof of \ref{_ergodicity_Theorem_}:}

The orbit $\Diff^+ \cdot \omega$ is dense in $\cF$ (in the torus
case) or, respectively, in $\cF^0$ (in the hyperk\"ahler case) if
and only if the orbit of the image of $\omega$ in $\Teich_s$ under
the action of
$\bGamma$ is dense in $\Teich_s$ (in the torus
case) or, respectively, in $\Teich_s^0$ (in the hyperk\"ahler case).

We claim that the latter condition holds if and only if the orbit $\bGamma\cdot
[\omega]$ is dense in $P$.

Indeed, in the hyperk\"ahler case this is true, since
$\Per: \Teich_s^0\to P$ is a
diffeomorphism. In order to prove the claim in the torus case, consider the subgroup $\Diff_H\subset\Diff^+$, consisting of all diffeomorphisms of $T^{2n}$ acting by identity on $H^* (T^{2n})$. The action of $\Diff_H$ on $\Teich_s$ preserves the fibers of $\Per$ and moreover is transitive on each fiber\footnote{Let us note that for $T^{2n}$, $n\geq 3$, the group $\Diff_H/\Diff_0$ is infinite (see \cite[Thm. 4.1]{_Hatcher_}, \cite[Thm. 2.5]{_Hsiang-Sharpe_}).} -- this follows from \ref{_Any_complex_structure_of_Kahler_type_on_torus_linear_Proposition_}
and the fact that two linear cohomologous forms on $T^{2n}=\R^{2n}/\Z^{2n}$ coincide. The latter transitivity yields the claim.

By \ref{_Ratner_for_H^2_Theorem_},
$\bGamma\cdot [\omega]$ is dense in $P$ if and only if the cohomology class $[\omega]$ is
not proportional to a rational one. We conclude that the orbit $\Diff^+ \cdot \omega$ is
dense in $\cF$ (in the torus case) or, respectively, in $\cF^0$ (in
the hyperk\"ahler case) if and only if the cohomology class $[\omega]$ is not proportional to a rational one.
This finishes the proof.
\endproof


\hfill


\noindent
{\bf Proof of \ref{_Ratner_for_H^2_Theorem_}:}

Both in the torus and in the hyperk\"ahler case set $k:=b_2:= \dim H^2 (M,\R)$.

In the case of the torus we identify $H^2(T^{2n}, \R)$ with the space $\big(\bigwedge^2
\R^{2n}\big)^*$ of linear 2-forms on $\R^{2n}$ and the action of
$\bGamma=\Diff^+/\Diff_0$ on $H^2(T^{2n}, \R)$ with
the action of $\Gamma :=SL(2n, \Z)$ on $\big(\bigwedge^2
\R^{2n}\big)^*$.
The latter action extends to the natural $SL (2n,\R)$-action on
$\big(\bigwedge^2
\R^{2n}\big)^*$. We also fix a linear identification of $H^2(T^{2n}, \R)=\big(\bigwedge^2
\R^{2n}\big)^*$ with
$\R^k$, so that $H^2 (T^{2n},\Q)\subset H^2 (T^{2n},\R)$ is identified with $\Q^k\subset \R^k$.

The stabilizer of a point $\eta\in P\subset H^2 (T^{2n},\R)$ under the
$SL (2n,\R)$-action is a closed Lie subgroup $H_\eta \subset SL
(2n,\R)$ isomorphic to $Sp (2n,\R)\subset SL (2n,\R)$ by an inner
automorphism of $SL (2n,\R)$. (Indeed, $\Per$ identifies elements of
$P$ with linear symplectic forms on $\R^{2n}$ defining the same
volume form and the $SL (2n,\R)$-action on $\big(\bigwedge^2
\R^{2n}\big)^*$
clearly preserves the set of such forms).

In the hyperk\"ahler case we identify $H^2 (M,\R)$ with $\R^k$ so that
the Bogomolov-Beauville-Fujiki form $q$ defined in Section~\ref{_BBF_form_Section_}
is identified with the standard $(3,k-3)$-quadratic form on $\R^k$ (and $H^2 (M,\Q)\subset H^2 (M,\R)$ is identified with $\Q^k\subset\R^k$).
The group of linear automorphisms of
$H^2 (M,\R)$ preserving
\[ P=\Theta_h=\left\{\ \eta \in H^2(M, \R)\ \ |\ \
q(\eta, \eta)>0,\ \int_M \eta^n = 1\ \right\}
\]
is then identified with
$$SO(3, k-3)\cong O(H^2(M, \R), q) \cap
SL (H^2(M, \R),\R).$$
 The stabilizer of a point $\eta\in P$
under the action of $SO(3, k-3)$ is a closed Lie subgroup

$H'_\eta\subset SO(3, k-3)$ isomorphic to $SO(2, k-3)\subset
SO(3, k-3)$ by an inner automorphism of $SO(3, k-3)$.
Let $H_\eta\subset SO^+ (3,k-3)$ be the identity component of $H'_\eta$.
It is isomorphic to $SO^+ (2, k-3)$ by the same inner automorphism of $SO(3, k-3)$.

Since $\bGamma\cdot P = P$, there is a natural homomorphism
$\bGamma\to SO(3, k-3)$. As it was shown in \cite{_V:Torelli_},
the image of this homomorphism is an arithmetic lattice in $SO(3, k-3)$. The intersection of the latter lattice with
$SO^+(3, k-3)$ will be denoted by
$\Gamma\subset SO^+ (3, k-3)$.

We are going to apply \ref{_Ratner's_arithmetic_Claim_} to $G=SL(2n,
\R)\supset H:=H_\eta \cong Sp(2n, \R)$, $\Gamma = SL (2n,\Z)\subset
G$ in the torus case and to $G=
SO^+ (3, k-3)\supset H:=H_\eta \cong
SO^+ (2,k-3)$ and the arithmetic lattice $\Gamma\subset G$ in the
hyperk\"ahler case. Indeed, in both cases
$G$ is the identity component of a linear
algebraic $\Q$-group $G'$ ($G'=G$ in the torus case and $G'= SO (p,q)$ in the
hyperk\"ahler case)
that does not admit non-trivial $\Q$-characters
(since it has no normal subgroups except for the center which is
discrete -- see e.g. \cite{_Ragozin_PAMS72_});
$\Gamma$ is an arithmetic lattice (by its definition and \ref{_lattice-identity-component_Proposition_}). The group $H_\eta$ is a closed Lie subgroup of $G$ generated by
algebraic one-parametric unipotent subgroups: indeed, it is obtained by an inner automorphism of $G$ from
$Sp(2n, \R)$ in the torus case and from $SO^+(2, k-3)$ in the
hyperk\"ahler case; the latter groups are generated by their algebraic one-parametric unipotent subgroups, by
\ref{_SL-n-R-SO-p-q-generated-by-unipotents_Proposition_}, hence so is $H_\eta$.

Hence
\ref{_Ratner's_arithmetic_Claim_} can be applied.

In both cases $G\cdot P = P$ and it is not hard to show that $G$
acts transitively on $P$. Thus the orbit $G\cdot \eta = P$ is
identified with the homogeneous space $G/H_\eta$. Let $x_\eta :=
eH_\eta \in  G/H_\eta = P$, where $e\in G$ is the identity element.
The orbit $\Gamma\cdot x_\eta$ in $G/H_\eta = P$ is exactly the
orbit $\Gamma\cdot x \subset P$. Thus, in order to prove
\ref{_ergodicity_Theorem_} we need to show that
$\overline{\Gamma\cdot x_\eta}= G/H_\eta$ if and only if the cohomology class $\eta$ is
not proportional to a rational one.

We will need the following lemmas.


\hfill


\lemma
\label{_no-intermediate-subgroups-between-H-and-G_Lemma_}

There are no intermediate connected Lie subgroups $Sp (2n,\R)\subsetneq S\subsetneq SL (2n,\R)$ and
$SO^+ (p-1,q)\subsetneq S\subsetneq SO^+ (p,q)$
(for any $p,q\in\Z_{>0}$).


\hfill


\lemma
\label{_rationality-of-the-stabilizer_Lemma_}

The group $H_\eta$ is
the identity component of a $\Q$-subgroup
if and only if $\eta\in H^2(M,\R)$ is proportional to a rational
cohomology class.


\hfill


Postponing the proof of the lemmas let us finish the proof of the theorem.

It follows from \ref{_Ratner's_arithmetic_Claim_} that
\[
\overline{\Gamma\cdot x_\eta} = \Gamma S_\eta \cdot x_\eta = \Gamma S_\eta e H_\eta = \Gamma S_\eta H_\eta,
\]
where
$S_\eta$ is the
identity component of the
smallest $\Q$-subgroup of
$G'$
containing
$H_\eta$.

It follows from
\ref{_no-intermediate-subgroups-between-H-and-G_Lemma_}
that for any $\eta \in P$ there are exactly two possibilities: either $S_\eta = H_\eta$ or $S_\eta = G$.

In the first case, when $S_\eta = H_\eta$, we have
\[
\overline{\Gamma\cdot x_\eta}=\Gamma S_\eta H_\eta\subsetneq G.
\]
Indeed,
$\Gamma$ is a discrete (hence countable) subgroup of $G$,
$S_\eta H_\eta$ is
a proper closed Lie subgroup of $G$ (hence a set of measure $0$ in $G$).
Hence, the set $\Gamma S_\eta H_\eta$ is of measure $0$ and, in particular,
a proper subset of $G$.

In the second case, when $S_\eta = G$, we clearly have
\[
\overline{\Gamma\cdot x_\eta} = \Gamma S_\eta H_\eta = \Gamma G H_\eta = G.
\]

The case $S_\eta = H_\eta$  happens if and only if $\eta$ is proportional to a rational
cohomology class. Indeed, if $\eta$ is proportional to a rational class, then $H_\eta$ is clearly a $\Q$-subgroup of $G$ and therefore
$S_\eta = H_\eta$. Conversely, if $S_\eta = H_\eta$, then $\eta$ is proportional to a rational
cohomology class, by \ref{_rationality-of-the-stabilizer_Lemma_}.

This finishes
the proof of  \ref{_ergodicity_Theorem_}.
\endproof

\hfill

\noindent {\bf Proof\footnote{This proof is due to M.Gorelik -- we
thank her for communicating it to us.} of
\ref{_no-intermediate-subgroups-between-H-and-G_Lemma_}.}

Let us
denote by $\gog$ the Lie algebra of $SL (2n,\R)$ (respectively, $SO
(p,q)$) and by $\goh$ the Lie algebra of $Sp (2n,\R)$ (respectively,
$SO (p-1,q)$). In both cases $\goh$ is a Lie subalgebra of $\gog$,
and since we want to rule out a connected intermediate Lie subgroup,  it suffices to show that there is no intermediate Lie subalgebra
$\gos$ such that $\goh\subsetneq \gos \subsetneq \gog$.

Indeed, assume $\goh\subset \gos \subset \gog$ and let us show that either $\gos=\gog$ or $\gos=\goh$.
Denote by $\gogc := \gog\otimes \C$, $\gohc := \goh\otimes \C$, $\gosc := \gos\otimes \C$ the corresponding complex Lie algebras.
It is enough to show that either $\gosc=\gogc$ or $\gosc=\gohc$.

View $\gogc$, $\gohc$ and $\gosc$ as modules over $\gohc$ with
respect to the adjoint action. Then $\gosc/\gohc\subset \gogc/\gohc$
is an inclusion of $\gohc$-modules. Since $\gohc$ is semi-simple
(both in the toric and the hyperk\"ahler cases), it is enough to
show that $\gogc/\gohc$ is an irreducible $\gohc$-module.

\bigskip
\noindent
{\bf The case of $\gog = \mathfrak{so}(p,q)$.}

The complexification of $\gog=\mathfrak{so}(p,q)$ is $\gogc=\mathfrak{so}(p+q, \C)$.

The algebra $\mathfrak{so}(n+1, \C)$, $n\in \N$, is the algebra of complex skew-symmet\-ric $(n+1)\times (n+1)$-matrices.
A skew-symmetric $(n+1)\times (n+1)$-matrix is of the form
$$T_{a,v}:=
\left(\begin{array}{c|c}
a & v\\
\hline
-v^t & 0
\end{array}\right),$$
where $a$ is a skew-symmetric $n\times n$-matrix
and $v$ is a $1\times n$-column.

The embedding $\mathfrak{so}(n,\C)\subset \mathfrak{so}(n+1,\C)$
corresponds to $v=0$: $\mathfrak{so}(n,\C)=\{T_{a,0}\}$. Clearly, there is an isomorphism of vector spaces:
$$\gogc=\gohc\oplus E,\ \text{ where }E:=\{T_{0,v}\}\subset \gogc.$$
Let us show that $E\cong \gogc/\gohc$ is an irreducible $\mathfrak{so}(n,\C)$-module.
Indeed, any element of $\mathfrak{so}(n,\C)$ is $T_{a,0}$ and one readily sees that
$$[T_{a,0}, T_{0,v}]=T_{0,av}.$$
In particular, $(ad\, g) e\in E$ for each $g\in \mathfrak{so}(n,\C), e\in E$,
so $E$ is an $\mathfrak{so}(n,\C)$-submodule of $\gogc$. Moreover, $E$ is isomorphic
to the standard module (that is, the $n$-dimensional complex vector
space with the natural $\mathfrak{so}(n,\C)$-action: $a\cdot v:=av$ for $a\in \mathfrak{so}(n,\C)$).
One easily checks that the standard module is irreducible and therefore
$\gogc/\gohc$ is irreducible as required.

\bigskip
\noindent
{\bf The case of $\gog = \mathfrak{sl}(2n,\R)$.}

First, let us recall the following basic facts concerning representations of Lie algebras.
Let $\gok$ be a complex Lie algebra and let $V, W$ be  $\gok$-modules.

We define the adjoint action of $\gok$ on $Hom(V,W)$ as follows:
\[
\bigl((ad\, g)(\psi)\bigr)(v):=
g(\psi(v))-\psi(gv)\ \text{ where } g\in\gok, \psi\in Hom(V,W), v\in V.
\]
(This action is called adjoint since $(ad\, g)\psi:=[g,\psi]=g\psi-\psi g$).

We define the $\gok$-module structure on $V^*:=End(V,\C)$ by viewing the base field
$\C$ as the trivial $\gok$-module (that is, $gv=0$ for all $g\in \gok$, $v\in \C$). Then
the $\gok$-module structure is given by
$$(gf)(v)=-f(gv),\ \ \text{ where } g\in\gok, f\in V^*, v\in V.$$

If $V$ is finite-dimensional, then the natural isomorphism
$$V^*\otimes V\cong End(V)$$
of vector spaces, given by $(f\otimes v)(v'):=f(v')v$, is an isomorphism of $\gok$-modules.

If $V$ is finite-dimensional and admits a non-degenerate $\gok$-invariant bilinear form  $(\cdot,\cdot)$ (that is,
such that $(gv,v')+(v,gv')=0$ for all $g\in\gok$, $v,v'\in V$), then the canonical isomorphism
$\iota: V\cong V^*$, given by $\iota(V)(v')=(v,v')$, is a $\gok$-isomorphism.

Let us now return to the case $\gog = \mathfrak{sl}(2n,\R)$.
The complex Lie algebra $\gohc=\mathfrak{sp}(2n,\C)$ is of type $C_n$.
Let $R$ be the natural representation of $\mathfrak{sp}(2n,\C)$: it is a $2n$-dimensional complex vector
space with the natural action of $\mathfrak{sp}(2n,\C)$. One has $R^*\cong R$, since $R$ admits a non-degenerate invariant bilinear form
(the symplectic form).

Observe that $\gogc:=\fsl(2n,\C)$ is a $\gohc$-submodule in $End(R)=\fgl(2n,\C)=\gogc\oplus \C$,
where $\C$ is the trivial $\gohc$-module.
Thus, in order to show that $\gogc/\gohc$ is irreducible, it is enough to verify
that $End(R)$ is the sum of three irreducible representations
(they are automatically $\gohc$, $\gogc/\gohc$ and $\C$).

Indeed, we have
$$End(R)=R\otimes R=S^2(R)\oplus \Lambda^2(R),$$
where $S^2(R)$, $\Lambda^2(R)$ are, respectively, the symmetric and the exterior squares of $R$.
Using \cite[Table 5]{_Onischik-Vinberg_}, we obtain
$$S^2(R)=R(2\pi_1),\ \Lambda^2(R)=R(\pi_2)\oplus R(\pi_0),$$
where $R(2\pi_1), R(\pi_2), R(\pi_0)$ are some irreducible representations of $\mathfrak{sp}(2n,\C)$ (in fact, $R(\pi_1)=R$,
$R(\pi_0)=\C$ and $R(2\pi_1)=\gohc = \mathfrak{sp}(2n,\C)$ is the adjoint representation of $\mathfrak{sp}(2n,\C)$).

Thus, $End(R)= \gogc\oplus \C$ is a sum of three irreducible $\gohc$-modules and therefore
$\gogc$ is a sum of two irreducible $\gohc$-modules one of which is $\gohc$. Therefore $\gogc/\gohc$
is an irreducible $\gohc$-module as required. This finishes the proof of \ref{_no-intermediate-subgroups-between-H-and-G_Lemma_}.
\endproof


\hfill


\noindent
{\bf Proof of \ref{_rationality-of-the-stabilizer_Lemma_}.}

If $\eta$ is proportional to a rational class, then the stabilizer $H_\eta$ of $\eta$ is clearly a $\Q$-subgroup of $G$ and therefore is the identity
component of
a $\Q$-subgroup of $G$ (because $H_\eta$ is connected).
This proves the ``if" part.

Let us prove the ``only if" part.
Assume that $H_\eta$ is the identity component of a $\Q$-subgroup $S$.
Then the Lie algebras $\Lie (S)$ and $\Lie (H_\eta)$ of $S$ and $H_\eta$ coincide: $\Lie (S) = \Lie (H_\eta)$. Since $S$ is the common zero level set of a finite number of polynomials (in the coefficients of matrices in $G$) with rational coefficients, its Lie algebra $\Lie (S)$ (that is, the tangent space of $S$ at the identity) is the common zero level set of a finite number of linear functions with rational coefficients, hence a Lie subalgebra of the Lie algebra $\Lie (G)$ of $G$ admitting a basis formed by matrices with rational coefficients.

The action of the Lie group $G$ on $H^2 (M;\R)$ induces a Lie algebra action of $\Lie (G)$ on $\R^k=H^2 (M;\R)$, that is, a Lie algebra homomorphism $F:\Lie (G)\to \gogl (k,\R)$. Since $H_\eta$ is the stabilizer of $\eta\in \R^k=H^2 (M;\R)$ under the $G$-action on $\R^k$, the Lie algebra $F(\Lie (H_\eta)) = F(\Lie (S))$ is the Lie subalgebra of $\gogl (k,\R)$ annihilating the space $\Span_\R \{\eta\}$.
Since $\Lie (S)$ has a basis formed by matrices with rational coefficients, so does $F(\Lie (H_\eta)) = F(\Lie (S))$, because $F$ is defined over $\Q$.
Hence, $\Span_\R \{\eta\}$ has to be spanned over $\R$ by a rational vector. Thus, $\eta$ is proportional to a rational vector, that is, to a rational class in $H^2 (M;\R)$.

This finishes the proof of the ``only if" part and of the lemma.
\endproof


\hfill


\noindent {\bf Acknowledgements:} We are grateful to L.Polterovich
for useful discussions and many interesting suggestions. We thank
M.Gorelik for helping us out with
\ref{_no-intermediate-subgroups-between-H-and-G_Lemma_},
E.Opshtein
for a stimulating discussion on packings by polydisks, A.Nevo for a
clarification concerning \ref{_Ratner's_arithmetic_Claim_} and
D.McDuff, F.Schlenk and an anonymous referee for pointing out
inaccuracies in the first version of the paper. The present work
is part of the first-named author's activities within CAST -- a
research network program of the European Science Foundation (ESF).


\hfill


{\small
\noindent {\sc Michael Entov\\
Department of Mathematics,\\
Technion - Israel Institute of Technology,\\
Haifa 32000, Israel}\\
{\tt  entov@math.technion.ac.il}
}
\\

{\small
\noindent {\sc Misha Verbitsky\\
Laboratory of Algebraic Geometry,\\
National Research University HSE,\\
Faculty of Mathematics, 7 Vavilova Str. Moscow, Russia}\\
{\tt  verbit@mccme.ru}
}


{\small
\begin{thebibliography}{AV1}

\bibitem[AV]{_Americ-Verbitsky:Teich_s_}
Amerik, E., Verbitsky, M., {\em Teichmuller space for hyperkahler
and symplectic structures}, J.Geom.Phys \textbf{97} (2015), 44-50.



\bibitem[Bea]{_Beauville_}
 Beauville, A. {\em
Varietes K\"ahleriennes dont la premi\`ere classe de Chern est
nulle.}  J. Diff. Geom. \textbf{18} (1983), 755-782.

\bibitem[Bes]{_Besse:Einst_Manifo_}
Besse, A., {\em Einstein Manifolds}, Springer-Verlag, New York,
1987.



\bibitem[Bo1]{_Bogomolov:decompo_}
Bogomolov, F. A., {\em On the decomposition of K\"ahler manifolds
with trivial canonical class}, Math. USSR-Sb. \textbf{22} (1974),
580-583.

\bibitem[Bo2]{_Bogomolov:defo_}
F. Bogomolov, {\em Hamiltonian K\"ahler manifolds}, Sov. Math. Dokl.
\textbf{19} (1978), 1462-1465.

\bibitem[BH2]{_Borel-HC1962_}
Borel, A., Harish-Chandra, {\em Arithmetic subgroups of algebraic groups}, Ann. of Math. (2) \textbf{75} (1962), 485-535.

\bibitem[Cal]{_Calabi_} Calabi, E.,
{\em On K\"ahler manifolds with vanishing canonical class}, in {\em
Algebraic geometry and topology. A symposium in honor of S.
Lefschetz, 78-89}, Princeton University Press, Princeton, N. J.,
1957.

\bibitem[Cam]{_Campana:isotrivial_}
Campana, F., {\em Isotrivialit\'e de certaines familles
K\"ahl\'eriennes de vari\'et\'es non projectives}, Math. Z.
\textbf{252} (2006), 147-156.



\bibitem[CDV]{_CDV:threefolds_}
Campana, F., Demailly, J.-P., Verbitsky, M., {\em Compact K\"ahler
3-manifolds without non-trivial
  subvarieties},
 Algebr. Geom. \textbf{1} (2014), 131-139.


\bibitem[Cat]{_Catanese:moduli_}
Catanese, F., {\em A Superficial Working Guide to Deformations and
Moduli}, in {\em Handbook of moduli, Vol. I, 161-215}, Int. Press,
Somerville, MA, 2013.



\bibitem[DP]{_Demailly_Paun_}
Demailly, J.-P., Paun, M., {\em Numerical characterization of the
K\"ahler cone of a compact K\"ahler manifold}, Ann. of Math. \textbf{159} (2004), 1247-1274.

\bibitem[Don]{_Donaldson:ellipt_}
Donaldson, S. K.,
{\em Two-forms on four-manifolds and elliptic equations}, in {\em
Inspired by S. S. Chern, 153-172},
World Sci. Publ., Hackensack, NJ, 2006.


\bibitem[Dou1]{_Douady_1966} Douady, A.,
{\em Le probl\`eme des modules pour les sous-espaces analytiques compacts d'un espace analytique donn\'e}, Ann. de l'inst. Fourier, \textbf{16} (1966), 1-95.

\bibitem[Dou2]{_Douady_}
Douady, A., {\em Le probl\`eme des modules pour les vari\'et\'es analytiques complexes (d'apr\`es Masatake Kuranishi)},
in {\em S\'eminaire Bourbaki, Vol. 9, Exp. No. 277, 7-13}, Soc. Math. France, Paris, 1995.

\bibitem[FH]{_Fricke_Habermann:symp-moduli_}
Fricke, J., Habermann, L.,
{\em On the geometry of moduli spaces of symplectic structures},
     Manuscripta Math. \textbf{109} (2002), 405-417.



\bibitem[F1]{_Fujiki:compactness_}
Fujiki, A., {\em
Closedness of the Douady spaces of compact K\"ahler
spaces}, Publ. Res. Inst. Math. Sci. \textbf{14} (1978/79), 1-52.

\bibitem[F2]{_Fujiki:Douady_}
Fujiki, A., {\em Countability of the Douady space of a complex
space}, Japanese J. of Math. \textbf{5} (1079), 431-447.


\bibitem[F3]{_Fujiki:HK_}
Fujiki, A., {\em On the de Rham Cohomology Group of a Compact
K\"ahler Symplectic Manifold}, Adv. Stud.
Pure Math. \textbf{10} (1987), 105-165.


\bibitem[Gro]{_Gromov_} Gromov, M., {\it Pseudoholomorphic curves in
symplectic manifolds}, Invent. Math. \textbf{82} (1985), 307-347.





\bibitem[Ham]{_Hamilton:Nash_}
Hamilton, R.S., {\em The inverse function theorem of Nash and Moser},
Bull. AMS
\textbf{7} (1982), 65-222.

\bibitem[Hat]{_Hatcher_} Hatcher, A.E., {\em Concordance spaces, higher simple homotopy theory, and applications}, in {\em Proc. Sympos. Pure Math., \textbf{32}, 3-21}, AMS, Providence, RI, 1978.

\bibitem[HS]{_Hsiang-Sharpe_} Hsiang, W.C., Sharpe, R.W., {\em Parametrized surgery and isotopy}, Pacific J. Math. \textbf{67}
(1976), 401-459.


\bibitem[H]{_Huybrechts:lec_}
 Huybrechts, D.,
{\em Compact hyperk\"ahler manifolds},
in {\em Calabi-Yau manifolds and related geometries (Nordfjordeid, 2001), 161-225},
Springer-Verlag, Berlin, 2003.


\bibitem[Ko]{_Kodaira_}
Kodaira, K., {\em Complex manifolds and deformation of complex structures}, Springer-Verlag, New York, 1986.

\bibitem[KoSp]{_Kod-Spen-AnnMath-1960_}
Kodaira, K., Spencer, D.C., {\em On deformations of complex analytic structures. III. Stability theorems for complex structures},
Ann. of Math. \textbf{71} (1960), 43-76.


\bibitem[KSS]{_Kleinbock_etc:Handbook_}
Kleinbock, D., Shah, N., Starkov, A.,
{\em Dynamics of subgroup actions on homogeneous spaces of
  Lie groups and applications to number theory}, in {\em Handbook
of dynamical systems, Vol. 1A, 813-930}, North-Holland,
Amsterdam, 2002.


\bibitem[Ku]{_Kuranishi_}
Kuranishi, M., {\em On the locally complete families of complex analytic structures}, Ann. of Math. \textbf{75} (1962), 536-577.

\bibitem[LMcDS]{_LMcDS_}
Latschev, J., McDuff, D., Schlenk, F., {\em The Gromov width of
4-dimen\-sional tori}, Geom. and Topol. \textbf{17} (2013), 2813-2853.

\bibitem[McDP]{_McD-Polt_}
McDuff, D., Polterovich, L., {\em Symplectic packings and algebraic geometry.
With an appendix by Yael Karshon}, Invent. Math. \textbf{115} (1994), 405-434.


\bibitem[McDS]{_McD-Sal-intro_} McDuff, D., Salamon, D.,
{\em Introduction to symplectic topology},
2nd edition, Oxford Univ. Press, Oxford, 1998.



\bibitem[Mor]{_Morris:Ratner_}
Morris, D.W., {\em Ratner's Theorems on Unipotent Flows}, Univ. of
Chicago Press, Chicago, 2005.





\bibitem[Mos]{_Moser_} Moser, J., {\em On the volume elements on a
manifold}, Trans. AMS \textbf{120} (1965), 288-294.


\bibitem[OG1]{_OGrady1_} O'Grady, K.G.,
{\em Desingularized moduli spaces of sheaves on a K3}, J. Reine
Angew. Math. \textbf{512} (1999), 49-117.

\bibitem[OG2]{_OGrady2_} O'Grady, K.G.,
{\em A new six-dimensional irreducible symplectic variety}, J. of
Alg. Geom. \textbf{12} (2003), 435-505.

\bibitem[OV]{_Onischik-Vinberg_}
Onishchik, A.L., Vinberg, E.B.,  {\em Lie groups and algebraic
groups}, Springer-Verlag, Berlin, 1990.


\bibitem[Rag]{_Ragozin_PAMS72_} Ragozin, D.L.,
{\em A normal subgroup of a semisimple Lie group is closed}, Proc. of
AMS \textbf{32} (1972), 632-633.



\bibitem[Rat]{_Ratner_Duke_1991_} Ratner, M.,
{\em Raghunathan's topological conjecture and
  distributions
of unipotent flows}, Duke Math. J. \textbf{63} (1991),
235-280.

\bibitem[Rem]{_Remmert_MathAnn_1956_} Remmert, R.,
{\em Projektionen analytischer Mengen},
Math. Ann. \textbf{130} (1956), 410-441.

\bibitem[Sh]{_Shah:uniformly_}
Shah, N.A.,
{\em Uniformly distributed orbits of certain flows on homogeneous spaces,}
Math. Ann. \textbf{289} (1991), 315-333.

\bibitem[Sp]{_Spivak_}
Spivak, M., {\em Calculus on manifolds. A modern approach to
classical theorems of advanced calculus}. W. A. Benjamin, Inc., New
York-Amsterdam, 1965.

\bibitem[Ti]{_Tian_}
Tian, G.,
{\em Smoothness of the universal deformation space of compact Calabi-Yau manifolds and its Petersson-Weil metric}, in
{\em Mathematical aspects of string theory (San Diego, Calif., 1986), 629-646},
World Sci. Publishing, Singapore, 1987.


\bibitem[To]{_Todorov_} Todorov, A.N., {\em The Weil-Petersson geometry of the moduli space of $SU(n\geq 3)$ (Calabi-Yau) manifolds. I}, Comm. Math. Phys.
\textbf{126} (1989), 325-346.




\bibitem[V1]{_Verbitsky:Symplectic_II_}
Verbitsky, M., {\em Hyperk\"ahler embeddings and holomorphic
symplectic geometry II,}  Geom. and Funct. Analysis \textbf{5} (1995), 92-104.





\bibitem[V2]{_Verbitsky:Deforma_}
Verbitsky, M., {\em Deformations of
trianalytic subvarieties of
hyperk\"ahler manifolds}, Selecta Math. \textbf{4} (1998), 447-490.





\bibitem[V3]{_V:Torelli_}
Verbitsky, M.,
{\em A global Torelli theorem for hyperk\"ahler manifolds},
 Duke Math. J. \textbf{162} (2013), 2929-2986.

\bibitem[V4]{_Verbitsky:ergodic_}
Verbitsky, M., {\em Ergodic complex structures on hyperk\"ahler
manifolds}, pre\-print, arXiv:1306.1498, 2014.


\bibitem[V5]{_Verbitsky:ICM_}
Verbitsky, M.,
{\em Teichm\"uller spaces, ergodic theory and global Torelli theorem},
preprint, arXiv:1404.3847, to appear in
in {\em Proceedings of the International Congress of
Mathematicians (Seoul 2014), vol. II, 793-813}.





\bibitem[Voi]{_Voisin-Hodge_} Voisin, C., {\em Hodge theory and complex algebraic geometry I,II}.  Cambridge Univ. Press, Cambridge, 2002.




\bibitem[Wh]{_Whitney_}
Whitney, H., {\em Elementary structure of real algebraic varieties},
Ann. Math., \textbf{66} (1957), 545-556.




\bibitem[Yau]{_Yau_} Yau, S.T.,
{\em On the Ricci curvature of a compact K\"ahler manifold and the complex Monge-Amp\`ere equation. I},
Comm. Pure Appl. Math. \textbf{31} (1978), 339-411.


\end{thebibliography}
}
\end{document}